\documentclass[12pt]{article}
\usepackage{geometry,color}

\usepackage{verbatim}
\newif\ifdeveloping
\developingtrue

\ifdeveloping
\fi

\usepackage{amsfonts}
\usepackage{amssymb}
\usepackage{enumerate}
\usepackage{amsmath}
\usepackage{graphicx}

\newtheorem{ttt}{Theorem}[section]
\newtheorem{llll}[ttt]{Lemma}
\newtheorem{ccc}[ttt]{Proposition}
\newtheorem{eee}[ttt]{Example}
\newtheorem{fff}[ttt]{Fact}
\newtheorem{rrr}[ttt]{Remark}
\newtheorem{sss}[ttt]{Statement}
\newtheorem{ddd}[ttt]{Definition}
\newtheorem{qqq}[ttt]{Problem}
\newtheorem{cccc}[ttt]{Corollary}
\newtheorem{nnn}[ttt]{Notation}

\newcommand{\bt}{\begin{ttt}}
\newcommand{\bl}{\begin{llll}}
\newcommand{\bc}{\begin{ccc}}
\newcommand{\bex}{\begin{eee}}
\newcommand{\bfa}{\begin{fff}}
\newcommand{\br}{\begin{rrr}\upshape}
\newcommand{\bst}{\begin{sss}}
\newcommand{\bd}{\begin{ddd}\upshape}
\newcommand{\bq}{\begin{qqq}}
\newcommand{\bnn}{\begin{nnn}}
\newcommand{\bcor}{\begin{cccc}}

\newcommand{\bp}{\noindent\textbf{Proof. }}

\newcommand{\doboz}{\raisebox{3.5pt}{\fbox{{}}}}

\newcommand{\et}{\end{ttt}}
\newcommand{\el}{\end{llll}}
\newcommand{\ec}{\end{ccc}}
\newcommand{\eex}{\end{eee}}
\newcommand{\efa}{\end{fff}}
\newcommand{\er}{\end{rrr}}
\newcommand{\est}{\end{sss}}
\newcommand{\ed}{\end{ddd}}
\newcommand{\eq}{\end{qqq}}
\newcommand{\ecor}{\end{cccc}}
\newcommand{\enn}{\end{nnn}}

\newcommand{\ep}{\hspace{\stretch{1}}$\square$\medskip}

\newcommand{\lab}[1]{\label{#1}}

\newcommand{\al}{\alpha}

\newcommand{\om}{\omega}

\newcommand{\ka}{\kappa}

\newcommand{\iH}{\mathcal{H}}

\newcommand{\acal}{{\mathcal A}}

\newcommand{\hcal}{{\mathcal H}}
\newcommand{\ical}{{\mathcal I}}

\newcommand{\lcal}{{\mathcal L}}

\newcommand{\xcal}{{\mathcal X}}
\def\br#1;#2;{\left[ {#1} \right]^ {#2} }
\def\<{\left\langle}
\def\>{\right\rangle}
\newcommand{\setm}{\setminus}
\newcommand{\subs}{\subseteq}
\def\ss{\subseteq}
\newcommand{\empt}{\emptyset}
\newcommand{\Lim}{\mathrm{Lim}}

\newcommand{\conv}{\mathrm{conv}}
\newcommand{\tail}{\operatorname{tail}}

\newcommand{\zp}{z^+}
\newcommand{\zm}{z^-}

\newcommand{\scal}{{\mathcal S}}

\newcommand{\sm}{\setminus}

\newcommand{\es}{\emptyset}
\newcommand{\beq}{\begin{equation}}
\newcommand{\eeq}{\end{equation}}
\def\su{\subseteq}
\def\subs{\subseteq}

\def\<{\left\langle}
\def\>{\right\rangle}
\def\({\left(}
\def\){\right)}
\def\b{\left}
\def\j{\right}
\def\ve{\varepsilon}
\def\vt{\vartheta}
\def\sm{\setminus}
\def\real{\mathbb{R}}

\def\egesz{\mathbb{Z}}

\def\rac{\mathbb{Q}}
\def\rar{\rightarrow}
\def\ben{\begin{enumerate}}
\def\een{\end{enumerate}}
\def\Keq{\begin{equation}}
\def\Zeq{\end{equation}}
\def\bit{\begin{itemize}}
\def\eit{\end{itemize}}

\def\tp{\mathrm{tp}}
\def\cov{\mathrm{cov}}
\def\cf{\mathrm{cf}}

\def\ss{\subseteq}
\def\vp{\varphi}

\def\mc{\mathcal}
\def\es{\emptyset}

\def\gotO{\mathrm{On}}

\def\gotA{\mathfrak{A}}
\def\gotB{\mathfrak{B}}

\def\gotR{\mathfrak{R}}
\def\gotC{\mathrm{Card}}
\def\gotK{\mathbf{K}}
\def\goth{\mathbf{h}}
\def\gotJ{\mathbf{J}}
\def\gotI{\mathbf{I}}

\def\gotH{\mathbf{H}}

\def\prf{\textbf{Proof.}~}

\def\O{\mathrm{On}}

\newcommand{\dom}{{\rm dom}}

\newcommand{\stickT}{%
\setbox255=\hbox{\raise1ex\hbox{$\hspace{0.2pt}\,\bullet\,$}}
\mathord{\rlap{\hbox to\wd255{\hss\hbox{$|$}\hss}}
\box255}
}
\newcommand{\stickS}{%
\setbox255=\hbox{\raise0.6ex\hbox{$\scriptstyle\bullet$}}
\mathord{\rlap{\hbox to\wd255{\hss\hbox{$\scriptstyle|$}\hss}}
\box255}
}

\title{On splitting infinite-fold covers}

\author{M\'arton Elekes, Tam\'as M\'atrai, Lajos Soukup}

\begin{document}

\baselineskip=17pt

\maketitle

\begin{abstract}
Let $X$ be a set, $\ka$ be a cardinal number and let $\iH$ be a family of subsets of $X$
which covers each $x\in X$ at least $\ka$ times. What
assumptions can ensure that $\iH$ can be decomposed into $\kappa$ many disjoint subcovers?

We examine this problem under various assumptions on the set $X$ and on the cover $\iH$: among other situations, we consider covers of topological spaces by closed sets, interval covers of linearly ordered sets and covers of $\real^{n}$ by polyhedra and by arbitrary convex sets. We focus on these problems mainly for infinite $\kappa$. Besides numerous positive and negative results, many questions turn out to be
independent of the usual axioms of set theory.
\end{abstract}

\insert\footins{\footnotesize{MSC codes: Primary 03E05, 03E15
    Secondary 03C25, 03E04, 03E35, 03E40, 03E50, 03E65, 05C15, 06A05,
    52A20, 52B11}}
\insert\footins{\footnotesize{Key Words: splitting infinite cover, coloring,
    convex set, linearly ordered set, interval, closed set, compact set, cardinal,
    Cohen model, Continuum Hypothesis, Martin's Axiom}}

\section{Introduction}

Let $X$ be a set, $\ka$ and $ \lambda$ be cardinal numbers and let $\iH$ be a family of subsets of $X$
which covers each $x\in X$ at least $\ka$ times. What assumptions on $\iH$ can ensure that $\iH$ can be decomposed into $\lambda$ many disjoint subcovers? That is, which $\kappa$-fold cover can be split into $\lambda$ many subcovers?

Depending on personal taste, every mathematician can readily formulate the ``most relevant" context for the splitting problem; therefore splitting covers has a long-standing tradition. Unarguably, the most studied version of the problem is when $X$ is a topological space and $\iH$ is an open cover with special combinatorial properties. We do not attempt to summarize the vast amount of results  in this direction, the interested reader is referred to \cite{T} and the references therein. Nevertheless, we note that the literature on combinatorial properties of open covers is mainly concerned with how the combinatorics of open covers is related to topological properties of the underlying space. Therefore a  strong topological motivation for considering the given special classes of open covers is always present, and splitting is concerned as far as one is looking for ``nice" disjoint subcovers of a ``not so nice" open cover. Moreover, for most of the problems discussed in these papers the open covers are automatically countable.
In the present paper we do not work on open covers, and as we will see, we treat the problem of splitting covers from a more set theoretic point of view.

Another well-understood variant of the splitting problem deals with covers of finite  structures. The most interesting questions in this area ask for splitting the edge covers of (hyper)graphs, and almost optimal solutions of the relevant problems have  already been found long ago. But none of the available results concern infinite graphs or infinite-fold covers. In Section \ref{s:graph} we give a complete solution to the splitting problem of infinite-fold edge covers of graphs; we will recall the related finite combinatorial results there.

The situation turns out to be less clear if we are interested in the splitting of finite-fold covers of infinite sets, even in the seemingly simple case of covers of the plane by such familiar objects as circles, triangles or rectangles. To start with positive results, D.\ P\'alv\"olgyi and G.\ T\'oth  \cite{PaT} showed that for every open convex polygon $R$ in the plane there are constants $c(R)$ and $n(R)$ such that every $c(R)k^{n(R)}$-fold cover of the plane with translates of $R$ can be decomposed into $k$ disjoint covers (see also \cite{TT} for the special case of open triangles). The analogous decomposition result in the special case of centrally symmetric open convex polygons was obtained by J.\ Pach and   G.\ T\'oth \cite{PT}, and in \cite{ACC} it was shown that for such regions $n(R)=1$ can be chosen.   In contrast with these results, J.\ Pach, G.\ Tardos and   G.\ T\'oth \cite{PTT} constructed, for every $1 < k < \omega$, a $k$-fold cover of the plane (1) by open strips, (2) by axis-parallel open rectangles, (3) by the homothets of an arbitrary open concave quadrilateral which cannot be decomposed into two disjoint covers. In fact, D.\ P\'alv\"olgyi \cite{Pa} obtained a characterization of the open polygons for which a positive decomposition result holds.

However, the problem whether for a given convex subset $R$ of the plane there is a $k$ such that any $k$-fold cover of the plane with translates/homothets of $R$ can be decomposed into two disjoint subcovers is far from being solved. E.g.\ we do not know the answer when $R$ is an open or closed disk; but we remark that a positive answer may be hidden in a more than 100 page-long manuscript of J.\ Mani-Levitska and J.\ Pach. We also note that the situation in $\real^{3}$ can turn out to be completely different: in another unpublished work of J.\ Mani-Levitska and J.\ Pach, for every $k < \omega$ a $k$-fold cover of $\real^{3}$ with open unit balls is constructed which cannot be decomposed into two disjoint covers. D.\ P\'alv\"olgyi \cite{Pa} obtained the analogous negative answer for polyhedra in $\real^{3}$.

Our investigations were initiated by the question of J.\ Pach whether any infinite-fold cover of the plane by axis-parallel rectangles can be decomposed into two disjoint subcovers (see also \cite[Concluding remarks p.\ 12]{AHM}). After answering this question in the negative for $\omega$-fold covers, we started a systematic study of splitting infinite-fold covers in the spirit of J.\ Pach et al.; in the present paper we would like to publish our first results and state numerous open problems.

We have organized the paper to add structure as we go along. In Section \ref{s:arbit}, for any pair of cardinals $\kappa$ and $\lambda$, we study the splitting of covers of $\kappa$ by sets in $[\kappa]^{\leq \lambda}$. In Section \ref{s:graph}, we discuss the splitting of edge-covers of finite or infinite graphs. In the remaining sections of the paper we study covers by convex sets. In Section \ref{s:in}, we show that a cover of a linearly ordered set by convex sets is ``maximally" decomposable. After completing our work, it turned out that R.\ Aharoni, A.\ Hajnal and E.\ C.\ Milner \cite{AHM} obtained results earlier which are similar to our results in Section \ref{s:in}. Since our proofs are significantly simpler and yield slightly stronger results we decided not to leave them out.

In Section \ref{closedsets}, as a preliminary study to covers by convex sets on the plane,  we show that the splitting problem for covers by closed sets is independent of ZFC. Roughly speaking, under Martin's Axiom an indecomposable cover of $\real$ can be obtained even by the translates of one compact set; while in a Cohen extension of a model with GCH, every uncountable-fold cover by closed sets is ``maximally" decomposable. From these results, in Section \ref{s:convex} we easily get that the splitting problem for covers of $\real^{n}$ by convex sets is independent of ZFC. This independence is accompanied by two ZFC results.  We show that for very general classes of sets, including e.g.\ polyhedra, balls or arbitrary affine varieties, an uncountable-fold cover by such sets is ``maximally" decomposable.  On the other hand, we construct an $\omega$-fold cover of the plane by closed axis-parallel rectangles which cannot be decomposed into two disjoint subcovers. We close the paper with a collection of open problems.

\section{Terminology}

In this section we fix the notation which will be used in all of the forthcoming sections. We denote by $\gotO$ and $\gotC$  the class of ordinals and the class of cardinals, respectively. For any set $X$ and cardinal $\lambda$, $\mc{P}(X)$ denotes the power set of $X$, and $[X]^{\lambda}$, $[X]^{\leq \lambda}$ and $[X]^{< \lambda}$ stand for the families of those subsets of $X$
which have cardinalities $\lambda$, $\leq \lambda$ and $ < \lambda$, respectively. If $\kappa$ is an ordinal, $\Lim(\kappa)$ denotes the set of limit ordinals  $ < \kappa$. For $\alpha,\beta \in \gotO$ with $\alpha \leq \beta$, lett $[\alpha,\beta]=\{\gamma \in \gotO \colon \alpha \leq \gamma \leq \beta\}$. The ordinal interval $[\alpha,\beta)$ is defined analogously. If $f \colon X \rar Y$ and $A \ss X$ are given, $f[A] = \{f(x) \colon x \in A\}$.

When we consider covers of a set $X$, we do not want to exclude to use a set $H \ss X$ multiple times. This approach is motivated  both by theoretical and by practical reasons. First, the classical results for splitting finite-fold covers of finite graphs allow graphs with multiple edges, so it is reasonable to keep this generality while extending these results for infinite graphs and infinite-fold covers. Second, the natural operation of restricting a cover of $X$ to a subset of $X$ can easily result in a cover where some of the covering sets are used multiple times. Moreover, this generality does not cause any additional complication. The following definition makes our notion {\it cover} precise.

\bd\label{feed}
Let $X$ be an arbitrary set, let $\iH \ss \mc{P}(X)$ be an arbitrary family of subsets of $X$ and let $m \colon \mc{H} \rar \gotO \sm \{0\}$ be an arbitrary function. Then the \emph{cover on X by $\mc{H}$ with multiplicity $m$} is $\gotH = \{\{H\} \times m(H) \colon H \in \mc{H}\}$.
For $x \in X$, let $\gotH(x) = \{\langle H, \alpha \rangle \colon x \in H \in \iH , ~\alpha < m(H) \}.$ A cover is \emph{simple} if $m(H)=1$ $(H \in \iH)$; for simple covers we identify $\gotH$ with $\mc{H}$.

Let $Y \ss X$ and let $\ka$ be a cardinal number. Then $\gotH$ is a \emph{cover of $Y$} if
$|\gotH(x)| \ge 1$ for every $x \in Y$, and $\gotH$ is a \emph{$\ka$-fold cover of $Y$} if
$|\gotH(x)| \ge \ka$ for every $x \in Y$. We say 
$\gotH$ is a \emph{cover} (\emph{$\ka$-fold cover}, resp.) if it is a cover ($\ka$-fold cover, resp.) of $\cup \iH$.
\ed

In the sequel  $\iH$, $m$ and $\gotH$ will always be as in Definition \ref{feed}. To ease notation, the decomposition of a cover will be realized by coloring the covering sets.

\bd
\lab{d:good}
Let $\gotH$ be a cover on $X$, let $Y \ss X$ and let $\kappa \in \gotC$. A partial function $c : \gotH \to \gotO$ is a \emph{${\kappa}$-good coloring of $\gotH$ over $Y$}   if for every $x \in Y$ 
if $|\gotH(x)|\ge {\kappa}$ then ${\kappa} \ss c[\gotH(x)]$.
 Similarly, 
$c : \gotH \to \ka$ is a \emph{${\kappa}$-good coloring of $\gotH$}, or simply a \emph{${\kappa}$-good coloring}, if
 it is a ${\kappa}$-good coloring of $\gotH$ over $\cup \iH$.

The coloring  $c$ is a 
\emph{$[{\kappa},\infty)$-good coloring of $\gotH$ over $Y$}
if it is ${\lambda}$-good for each cardinal ${\lambda}\ge {\kappa}$. We say that $c$ is a \emph{maximally good coloring of $\gotH$ over $Y$}
if it is $1$-good and $[{\omega},\infty)$-good. The notions \emph{$[{\kappa},\infty)$-good coloring of $\gotH$} and \emph{maximally good coloring of $\gotH$} are defined analogously.
\ed

Clearly, a ${\kappa}$-fold cover of $Y$ has a $\ka$-good coloring over $Y$ if and only if
it can be partitioned into $\kappa$ many covers of $Y$.

Next we prove an easy  reduction theorem. 

\bc \label{mred} Let $X$ be a set, $Y \ss X$, $\kappa \in \gotC \sm \omega$, and let $\mc{H}$, $m$ 
and $\gotH$ be as in Definition \ref{feed}.
If there is a ${\kappa}$-good  coloring of $\mc{H}$ over $Y$ then there is a
  ${\kappa}$-good coloring of $\gotH$ over $Y$. The analogous result holds for $[{\kappa},\infty)$-good and maximally good colorings, as well.
\ec

 Proposition \ref{mred} allows us to consider only simple covers. We will
frequently use this reduction steps in inductive proofs: it allows
us to use the inductive assumption for multicovers  while 
we prove our statement in the special case of simple
covers. We will always state explicitly when
this reduction is used. Note that the assumption $\kappa$ being infinite cannot be left out: e.g.\ the complete graph on three vertices $G$ has an $n$-good edge coloring for every $n \geq 3$, while if we take each edge of $G$ with multiplicity two, then the resulting graph has no $4$-good edge coloring.

\medskip

\textbf{Proof of Proposition \ref{mred}.}
First suppose $c_{0} \colon \iH \rar \gotO$ is a ${\kappa}$-good coloring of
$\mc{H}$ over $Y$. Let $\chi \colon \O \sm \{0\} \rar \O$ satisfy $\chi[\lambda\sm\{0\}] = \lambda$ $(\lambda \in \gotC \sm \omega)$. For every $H \in \dom(c_{0})$ we define $$c(\<H,\alpha\>)=\b\{\begin{array}{ll}c_{0}(H) & \textrm{if } 
\alpha = 0, \\ \chi(\alpha) & \textrm{if } 0 < \alpha < m(H).
\end{array} \j.$$ Then for every $x \in Y$, $c_{0}[\iH(x)] \ss c[\gotH(x)]$; and for every $H \in \iH(x)$, $\omega \leq m(H)$ implies $m(H) \ss c[\gotH(x)]$. Hence $c$ is a ${\kappa}$-good coloring of
$\gotH$ over $Y$.
For $[{\kappa},\infty)$-good and maximally good colorings the proof  is identical.
\ep

In some inductive proofs we will also use the following lemma.

\bl\label{l:ha} Let $\kappa \in \gotC$ and $\alpha \in \O$ satisfy ${\omega}_1 \leq \kappa  = |{\alpha}|$. 
Then there is a function
$h_{\alpha}:{\kappa}\to {\alpha}$ such that
$h_{\alpha}[{\kappa}]={\alpha}$ and ${\kappa}' \ss h_{\alpha}[{\kappa}']$
for each  cardinal ${\omega}_1\le {\kappa}'<{\kappa}$.  
\el

\bp Let  $h_{\alpha}$ be such that $h_{\alpha}({\nu})=\tp(\Lim(\kappa)\cap {\nu})$ for ${\nu}\in
\Lim({\kappa})$
and $h_{\alpha}[{\kappa}\setm \Lim({\kappa})]={\alpha}$. It is obvious that $h_{\alpha}[{\kappa}]={\alpha}$. Let $\kappa' \in \gotC$ satisfy ${\omega}_1\le {\kappa}'<{\kappa}$.  Since for every $\beta \in [0,\kappa')$ there is a $\nu \in [0,\kappa')$ such that $\tp(\Lim(\kappa')  \cap \nu)=\beta$, we get ${\kappa}' \ss h_{\alpha}[{\kappa}']$. So
$h_{\alpha}$ satisfies the requirements.
\ep

Later on, we will apply Lemma \ref{l:ha} the following way. Suppose we have an $[{\omega}_1,\infty)$-good coloring $c:\gotH\to {\kappa}$ and an ordinal  
${\kappa}\le {\alpha}<{\kappa}^+$. Then the coloring 
$h_{\alpha}\circ c$ is an
$[{\omega}_1,\infty)$-good coloring
 with the additional property 
\begin{enumerate}[(m)]
\item if $|\gotH(x)|={\kappa}$ then ${\alpha} \ss (h_{\alpha}\circ
  c)[\gotH(x)]$.  
\end{enumerate}

\subsection{Hunting for the strongest possible decomposition result}

In this subsection we investigate how  one can obtain colorings 
with even stronger decomposition properties from
maximally good or $[{\omega}_1,\infty)$-good colorings.
The notions and results of this subsection are not used later
in the paper. The strongest possible decomposition notion is formulated in the
following terminology. 

\bd
\lab{d:maxim}
Let $\gotH$ be a cover of $X$ and let $Y \ss X$. Let $\goth \colon
\gotC \rar \gotO \sm \{0\}$ be a partial function satisfying
$\goth(\kappa) < \kappa^{+}$ ($\kappa \in  \dom (\goth))$. A partial
function $c : \gotH \to \gotO$ is an \emph{$\goth$-good coloring of
  $\gotH$ over $Y$} if for every $x \in Y$,
\begin{enumerate}[(m$^{+}$)]
\item \label{max+} if $|\gotH(x)| \in \dom(\goth)$ then
$\goth(|\gotH(x)|)\subseteq c[\gotH(x)]$.
\end{enumerate} \ed

\bc \label{crewd} Let $X$, $Y$, $\gotH$ and $\goth$ be as in
Definition \ref{d:maxim}. Suppose that 
\begin{enumerate}[(h1)]
\item \label{h1} if $\dom(\goth) \cap \omega$ is unbounded in $\omega$ then $\goth(\omega) = \omega$;
\item \label{h2} for each ${\mu} \in \gotC$ the
set 
$\{{\nu} \in \dom(\goth) \cap {\mu} \colon {\nu} < \goth({\nu})\}$ is not stationary in ${\mu}$.
\end{enumerate} Set $\mathbf{i}\colon \dom(\goth) \rar \gotO$, $\mathbf{i}(\kappa)=\kappa$ $(\kappa \in \dom(\goth))$.
If there exists an $\mathbf{i}$-good coloring of $\gotH$ over $Y$ 
 then there exists an $\goth$-good coloring of $\gotH$ over $Y$, as well.
\ec

Observe that a maximally good coloring is $\mathbf{i}$-good for $\dom(\mathbf{i}) = \{1\} \cup (\gotC \sm \omega)$. Before proving Proposition \ref{crewd}, we need a lemma in advance.
This lemma is a far fetching generalization of lemma \ref{l:ha}.

\bl\label{red} Let $\lambda$ be an infinite cardinal. Let the partial function
$\goth\colon \gotC \cap \lambda^{+} \rar \gotO$ satisfy (h\ref{h1}), (h\ref{h2}) and ${\kappa} \leq
\goth({\kappa}) < {\kappa}^{+}$ $({\kappa} \in  \dom(\goth))$. Then there exists a function $\chi \colon \lambda \rar \gotO$
such that 
 \begin{enumerate} 
 \item\label{w1} $\chi(n)=n$ for $n < 1+\sup (\dom(\goth) \cap \omega)$;
 \item\label{w2} $[0,\goth({\kappa})) \ss \chi[[0, {\kappa})]$  for ${\kappa} \in  \dom(\goth) \cap [\omega, \lambda]$. 
\end{enumerate}
\el
\prf By defining $\goth$ on $\gotC \cap [\omega, \lambda] \sm \dom(\goth)$ as $\mathbf{i}$, we can assume $\gotC \cap [\omega, \lambda] \ss \dom(\goth)$. We prove the statement by induction on $\lambda$. For $\lambda =
\omega$ a bijection $\chi \colon \omega \rar \goth(\omega)$ 
with $\chi(n)=n$ ($n < 1+\sup (\dom(\goth) \cap \omega)$) does the job. 

Let now $\lambda > \omega$ and suppose that the statement holds for
every cardinal ${\kappa} < \lambda$. If $\lambda$ is a successor,
say $\lambda = {\kappa}^{+}$, set $\goth_{{\kappa}} =
\goth|_{{\kappa}^{+}}$ and let $\chi_{{\kappa}}$ satisfy \ref{w1} and 
$[0,\goth_{{\kappa}}({\kappa}')) \ss \chi_{{\kappa}}[[0, {\kappa}')]$
    for every ${\kappa}' \in  \gotC \cap [\omega,
      {\kappa}]$. Define $\chi \colon \lambda \rar \gotO$ by
    $\chi|_{{\kappa}} = \chi_{{\kappa}}$, and $\chi|_{ [{\kappa},
	{\kappa}^{+})} \colon  [{\kappa}, {\kappa}^{+}) \rar
	[{\kappa}, \goth({\kappa}^{+}))$ being any bijection. Then
	  $\goth$ clearly fulfills the requirements. 

If $\lambda$ is a limit cardinal, take a strictly increasing
continuous cofinal sequence $\< \lambda_{\alpha} < \lambda \colon
\alpha < \cf(\lambda)\>$ of infinite cardinals such that
$\goth({\lambda}_{\alpha})={\lambda}_{\alpha}$ 
for every ${\alpha}<\cf({\lambda})$. Let $\goth(\lambda) \sm \lambda=
\bigcup\{ K_{\alpha} \colon \alpha < \cf(\lambda)\}$ such that for
every $\alpha \leq \alpha' < \cf(\lambda)$ we have $K_{\alpha} \ss
K_{\alpha'}$ and $|K_{\alpha}|\leq \lambda_{\alpha}$. For every
$\alpha < \cf(\lambda)$ define $\goth_{\alpha} \colon \gotC \cap
\lambda_{\alpha}^{+} \rar \gotO$, $\goth_{\alpha} = \goth|_{\gotC \cap
  \lambda_{\alpha}^{+}}$. By the inductive hypothesis, for every
$\alpha < \cf(\lambda)$ we have $\chi_{\alpha}\colon \lambda_{\alpha}
\rar \gotO$ 
 satisfying \ref{w1} and  $[0,\goth_{\alpha}({\kappa})) \ss \chi_{\alpha}[[0,
     {\kappa})]$ $({\kappa} \in  \gotC \cap [\omega,
       \lambda_{\alpha}])$. 

 Set $\chi|_{[0,\lambda_{0})} = \chi_{0}$; then \ref{w1} holds, and \ref{w2} holds for ${\kappa} \in  \gotC \cap [\omega, \lambda_{0}]$. For every $\alpha < \cf(\lambda)$ fixed, we define $\chi\colon
 [\lambda_{\alpha}, \lambda_{\alpha+1}) \rar \gotO$ as follows. Let
 $\vt \colon [\lambda_{\alpha},\lambda_{\alpha}\dot + |K_{\alpha}|)
 \rar K_{\alpha}$ be a bijection and let $\ve \colon [\lambda_{\alpha}
 \dot + |K_{\alpha}|, \lambda_{\alpha+1}) \rar \lambda_{\alpha+1}$ be
 the enumeration. Observe that  for every ${\kappa} \in \gotC \cap
 (\lambda_{\alpha}, \lambda_{\alpha+1})$ we have  $\ve({\kappa}) =
 {\kappa}$. 
 For $\eta \in [\lambda_{\alpha}, \lambda_{\alpha} \dot +
 |K_{\alpha}|)$ set $\chi(\eta) = \vt(\eta)$ while for $\eta \in
 [\lambda_{\alpha}\dot + |K_{\alpha}|, \lambda_{\alpha+1})$ set
 $\chi(\eta) = \chi_{\alpha+1}(\ve(\eta))$. Then for every ${\kappa}
 \in \gotC \cap (\lambda_{\alpha}, \lambda_{\alpha+1}]$ we have
 $\chi[[0, {\kappa})] \supseteq K_{\alpha}$ and $$\chi[[0, {\kappa})]
 \supseteq \chi_{\alpha+1}[\ve[[\lambda_{\alpha}\dot +
 |K_{\alpha}|,{\kappa})]] \supseteq \chi_{\alpha+1}[[0,{\kappa})]
 \supseteq [0,\goth_{\alpha+1}({\kappa})) =[0, \goth({\kappa})).$$ For
 every limit ordinal $\alpha < \cf(\lambda)$ we have $\lambda_{\alpha}
 = \sup_{\beta < \alpha}\lambda_{\beta}$, hence  $\lambda_{\alpha} \ss
 \chi[[0,\lambda_{\alpha})]$ and we have $\chi[[0, \lambda)] \supseteq
 [0, \goth(\lambda))$ as well. This completes the proof. 
 \ep

In Lemma \ref{red}, the technical assumptions (h\ref{h1}) and (h\ref{h2}) cannot be left out. This is obvious for (h\ref{h1}). For (h\ref{h2}), observe that if
$\lambda$ is a Mahlo cardinal then for the function $\goth(\nu) = \nu \dot +
1$ $(\nu \in \gotC \cap \lambda)$, where $\dot +$ denotes ordinal addition, a function $\chi$ satisfying
 conclusion \ref{w2} of Lemma \ref{red} would induce a regressive function on the set $\gotC \cap \lambda$, which is stationary in $\lambda$. This explains the assumptions on
$\goth$ in Proposition \ref{crewd}. 

\medskip

\textbf{Proof of Proposition \ref{crewd}.} Let $c_{\mathbf{i}} \colon
\gotH \rar \gotO$ be an $\mathbf{i}$-good coloring of $\gotH$ over
$Y$. Set $\lambda = \sup \dom (\goth)$, $$\goth^{+}(\kappa) =
\left\{\begin{array}{ll} \max\{\goth(\kappa), \kappa\}, & \kappa \in
\dom(\goth), \\ \kappa, & \kappa \in \lambda^{+} \sm \dom(\goth);
\end{array} \right.$$ and let $\chi \colon \lambda \rar \gotO$ be the
function of Lemma \ref{red} for $\goth^{+}$. Define $c\colon \gotH
\rar \gotO$ by $c= \chi \circ c_{\mathbf{i}}$. Then $c$ is clearly an
$\goth$-good coloring of $\gotH$ over $Y$. 
\ep

\section{Arbitrary sets}
\label{s:arbit}

In this section we briefly recall some easy coloring results which hold for covers by arbitrary sets. Let us recall that Axiom Stick is the statement 
\begin{multline}\notag
\textrm{
 there is a family } \scal\subs [{\omega}_{1}]^{{\omega}} \textrm{ such that } \\
|\scal|={\omega}_1 \textrm{ and } \forall X\in [{\omega}_1]^{{\omega}_1} ~\exists S\in \scal ~( S\subs X). 
  \end{multline}
We refer to \cite{Ku2} for the definition of Martin's Axiom for various posets.

\bt\label{asse} ~
  \begin{enumerate}
  \item\label{small} 
Let $\mu$ be an arbitrary cardinal and $\kappa$ be an infinite 
cardinal. Then every  $\kappa$-fold cover $\gotH$ with $\iH \subs [\mu] ^{\leq \kappa}$  
has a  ${\kappa}$-good coloring.

\item\label{power} For each  infinite cardinal ${\mu}$ 
there is a simple ${\mu}$-fold cover $\gotH\subs [ 2^{\mu}]^{2^{\mu}}$ of $2^{\mu}$
with $|\gotH|={\mu}$ which does not have a 2-good coloring.

\item\label{cohen} $MA_{{\mu}}(Fn({\mu},2;{\omega}))$ implies
that 
every   cover $\gotH $ with $\iH \subs [ {\mu}]^{\leq \mu}$
and  $|\gotH|\le {\mu}$ 
has an  ${\omega}$-good coloring.

\item\label{countable} $MA_{{\mu}}(Fn({\omega},2;{\omega}))$ implies
that 
every  cover $\gotH $ with $\iH \subs [ {\mu}]^{\leq \mu}$
and  $|\gotH|\leq {\omega}$   
has an  ${\omega}$-good coloring.

\item\label{stick}  Axiom Stick implies
that 
there is a simple
cover $\gotH\subs [{\omega}_1]^{{\omega}_1}$
with $|\gotH|={\omega}_1$ and 
$|\gotH({\xi})|={\omega}$ $({\xi}\in {\omega}_1 )$
which does not have a 2-good coloring.

\item\label{fss} It is consistent that $2^{\omega}$ is arbitrarily large and the following two statements hold:
  \begin{enumerate}[(i)]
   \item there is a simple  
cover $\gotH\subs [{\omega}_1]^{{\omega}_1}$
with $|\gotH|={\omega}_1$ and $|\gotH({\xi})|={\omega}$ 
$({\xi}\in {\omega}_1)$ 
which 
does not have a 2-good coloring;
\item
for each ${\mu}<2^{\omega}$
every  cover $\gotH $ with $\iH \subs [{\mu}]^{\leq {\mu}}$
and  $|\gotH|= {\omega}$  
has an  ${\omega}$-good coloring.

  \end{enumerate}
  \end{enumerate}
\et

\bp
\ref{small}. Set $X = \bigcup \mc{H}$. If ${\mu}\leq {\kappa}$, a straightforward transfinite induction yields the required coloring, as follows. Let $\vp \colon X \times \kappa \rar \kappa$ be a bijection. By induction, for every $\alpha < \kappa$ we define partial colorings $c_{\alpha} \colon \gotH \rar \kappa$ satisfying $|c_{\alpha}|=1$, as follows. If $\alpha < \kappa$ and $c_{\beta}$ is defined for $\beta < \alpha$, let $x \in X$, $\chi < \kappa $ satisfy $\vp(x,\chi) = \alpha$. Since  
  $|\gotH(x)|\geq {\kappa}$, we can pick an $\langle H , \gamma \rangle \in \gotH\sm \bigcup_{\beta < \alpha} \dom(c_{\beta})$ satisfying $x \in H$ and set $c_{\alpha}(\langle H , \gamma \rangle) = \chi$. This completes the $\alpha^{\textrm{th}}$ step of the construction. Then $c = \bigcup_{ \alpha < \kappa} c_{\alpha}$ is a $\kappa$-good coloring of
$\gotH$. 

Suppose now $\mu>{\kappa}$.
For each ${\xi}\in {\mu}$ let $
\gotH'_{\xi} \in [\gotH({\xi})]^{\kappa}$.
Let $\equiv$ be the equivalence relation on ${\mu}$ generated by the
 relation
 $\{({\xi},{\zeta}) \colon {\xi},{\zeta}\in {\mu}, ~\gotH'_{\xi}\cap
 \gotH'_{\zeta} \neq \es\}$.
Since $\iH \ss [{\mu}]^{\leq \kappa}$ and $|\gotH'_{\xi}|={\kappa}$
for each ${\xi}\in {\mu}$, the equivalence classes
 $\{X_i \colon i\in I\}$ of $\equiv$ have cardinalities $\le{\kappa}$.
Set $\gotH_i=\cup\{\gotH'_{\xi} \colon {\xi}\in X_i\}$ $(i  \in I)$. Then $|X_i| \leq |\gotH_{i}| =
{\kappa}$ and $|\gotH_{i}(x)|={\kappa}$ for each  $x\in X_{i}$. So by
the inductive hypothesis,
there is a ${\kappa}$-good coloring $c_i \colon \gotH_i \rar {\kappa}$
over $X_{i}$ $(i \in I)$. Then $c =\bigcup_{i\in I}c_i$ is a good
$\kappa$-coloring over $X$.

\ref{power}. Set $X = \{f\in 2^{\mu}:|f^{-1}\{1\}|={\mu}\}$, and for every $\alpha < \mu$ let
$H_{\alpha} = \{x \in X \colon x(\alpha) = 1\}$. Then $|X| = 2^{\mu}$ and $\gotH =
\{H_{\alpha} \colon \alpha < \mu\}$ is a simple cover of $X$
with  the required properties.

\ref{cohen}. Consider the poset 
$ \mc{P} = \{(\mc{K}, c) \colon \mc{K} \in [\gotH]^{<\omega}, ~c
\colon \mc{K} 
\rar \omega\}$ with partial order $(\mc{K}', c') \leq (\mc{K}, c) $ if
and only if $\mc{K} \ss \mc{K}'$ and $c \ss c'$. Observe that for
every $\alpha \in \mu$ with $|\gotH(\alpha)| \geq \omega$ and $\chi <
\omega$, 
the set 
$$D_{{\alpha},\chi}=\{(\mc{K}, c) \in \mc{P} \colon \exists \langle H , \gamma
\rangle \in \gotH ~(\alpha \in H \textrm{ and }c(\langle H , \gamma
\rangle) = \chi)\}$$ is dense in $\mc{P}$. So by applying
$MA_{{\mu}}(Fn({\mu},2;{\omega}))$ to the poset $\mc{P}$ and 
the family $\mc{D}=\{D_{{\alpha},{\lambda}} \colon {\alpha}\in
{\mu},~\chi<{\omega} \}$
of dense sets of cardinality ${\mu}$, the statement 
follows. The argument for \ref{countable} is similar.

\ref{stick}. Let $\mc{S} \ss [\omega_{1}]^{\omega}$ be a stick
family. Set $\mc{X} = \{S \cup \{\alpha\} \colon S \in \mc{S}, ~\alpha  < \omega_{1}\}$; then $\mc{X}$ is also a stick family with the additional properties that 
\ben[(i)]
\item\label{fr1} $|\mc{X}(\alpha)| = \omega_{1}$ $(\alpha < \omega_{1})$; 
\item\label{fr2} for every $\alpha, \beta < \omega_{1}$ with $\alpha \neq \beta$ there is an $X \in \mc{X}$ with $\alpha \in X$ and  $\beta\notin X$.
\een
For every $\alpha < \omega_{1}$ let $H_{\alpha} = \{X \in
\mc{X} \colon \alpha \in X\}$. Since $|\mc{X}| = \omega_{1}$, it is enough to show that $\gotH = \{H_{\alpha} \colon
\alpha < \omega_{1}\} $ is a cover of $\mc{X}$ with the required properties. By (\ref{fr1}), $|H_{\alpha}|= \omega_{1}$  $(\alpha < \omega_{1})$. By (\ref{fr2}), $H_{\alpha} \neq H_{\beta}$ for $\alpha \neq \beta$., i.e.\ $\gotH$ is a simple cover. For every $X\in \mc{X}$ we have $|\gotH(X)|=|\{H_{\alpha} \colon {\alpha}\in X\}|=|X|={\omega}$. Finally assume that for some $I \ss \omega_{1}$, 
$\gotH_0=\{H_{\alpha}\colon {\alpha}\in I\}$ covers $\mc{X}$. This means
 $I \cap X\ne \empt$ for each $X\in \mc{X}$.
Then $|{\omega}_1\setm I|\le {\omega}$ because every 
$J\in [ {\omega}_1]^{{\omega}_1}$ contains some element of
$\mc{X}$. Thus $X\subs I$ for some $X\in \mc{X}$ so this 
$X$ is not covered by $\gotH\setm \gotH_0$.

\ref{fss}.
In \cite{FSS} it was proved that it is consistent that  $2^{\omega}$ is arbitrarily large,
Axiom Stick holds and $MA(Fn({\omega},2;{\omega}))$ also holds. So by \ref{countable} and \ref{stick}, that model satisfies the requirements of \ref{fss}.
\ep

Theorem \ref{asse} is to be compared with the results of Section \ref{closedsets} on splitting closed covers.

\section{Graphs}
\lab{s:graph}

Now we investigate the interesting special case of graphs, that
is, when each covering set has $2$ elements. In this section ``graph" means an undirected possibly infinite graph where multiple edges are allowed, but we exclude loops. We follow the standard notation, i.e.\ $G = (V,E)$ denotes the graph with vertex set $V$ and edge set $E$. For every $V' \ss V$, $G[V']$ denotes the subgraph of $G$ spanned by $V'$. According to our convention, for $v,w \in V$, the set of edges containing $v$ is denoted by $E(v)$ and $E(v,w)$ stands for the edges connecting $v$ to $w$. For every $v \in V$, $d_{G}(v)$ stands for the degree of $v$ in $G$, i.e.\ $d_{G}(v) = |E(v)|$ where multiple edges are counted with multiplicity. Set $\Delta(G) = \sup\{d_{G}(v) \colon v \in V\}$; the supremum of edge multiplicities is denoted by $\mu(G)$. For every $E' \ss E$, $V[E']$ is the set of vertices of the edge set $E'$. A graph $G$ is \emph{$n$-regular} if $d_{G}(v) = n$ $(v \in V)$. A \emph{complete matching} in $G$ is a subgraph of $G'$ of $G$ satisfying $d_{G'}(v) = 1$ $(v \in V)$.

As we mentioned in the introduction, the splitting problem for finite graphs is much studied (see e.g.\ \cite[Chapter 28]{Sc}), and the following result,  originally due to R.\ P.\ Gupta \cite[Theorem 2.2 p.\ 500]{G}, solves our problem for finite graphs.

\bt\label{Gupta}
Let $1 \le n <\om$.
Let $G = (V,E)$ be a finite graph and let $X \su V$ be such that for every $x \in X$ we have $d_{G}(x)  \geq n+\mu(G)$. Then $E$ has an $n$-good coloring over $X$.
\et

The main result of this section is the extension of Theorem \ref{Gupta} for
infinite graphs such that, in addition, for the existence of 2-good colorings a
necessary and sufficient condition is given. First we show that even for simple graphs, in order to ensure the existence of
an $n$-good coloring, the condition on the degree of vertices cannot be weakened to
$d_{G}(x) \geq n$. To see this, we will use the following constructions proposed by Gyula Pap.

For every $n < \omega$ let $K_{n}$ denote the complete graph on the vertex set
$\{v_{0}, v_{1}, \dots, v_{n-1}\}$. For odd $n$ let $K^{-}_{n}$ denote the graph obtained
from $K_{n}$ by \emph{deleting} the edges $\{v_{0}, v_{n-1}\}$ and $\{v_{2k},v_{2k+1}\}$ $(k <
(n-1)/2)$. Take two disjoint copies of $K^{-}_{n+2}$, say on the vertex sets
$\{v_{0}, v_{1}, \dots, v_{n+1}\}$ and $\{v'_{0}, v'_{1}, \dots, v'_{n+1}\}$ and let $D_{n}$ denote
the graph obtained as the union of the two copies of $K^{-}_{n+2}$ and the edge
$\{v_{0}, v'_{0}\}$.

\bc\label{Pet} Let $2 \leq n < \omega$.

\ben
\item\label{Pet1} If $n$ is  even, $K_{n+1}$ is an $n$-regular graph with no $n$-good coloring.
\item\label{Pet2} If $n$ is odd, $D_{n}$ is an $n$-regular graph with no $n$-good coloring.
\een
\ec
\bp
We prove the statements simultaneously. It is obvious that $K_{n+1}$ and
$D_{n}$ are $n$-regular graphs. Hence an $n$-good coloring of any of these
graphs is a partition of their edge set into $n$ disjoint complete matchings.

Now for even $n$, $K_{n+1}$ has no complete matchings at all since its vertex set has
odd cardinality. Also for cardinality reasons, if $n$ is odd any complete matching of $D_{n}$ must
contain the edge $\{v_{0}, v'_{0}\}$. Hence $D_{n}$ has no two disjoint
complete matchings; in particular, by $2 \leq n$, the $n$-regular graph $D_{n}$ has no $n$-good coloring.
\ep

We also note that Theorem \ref{asse}.\ref{small}   completely solves the splitting
problem  for infinite-fold edge-covers.

\bt\label{UJU}
Let $G = (V,E)$ be a graph, let $X \ss V$ be arbitrary and suppose that for every $x \in X$ we have $d_{G}(x) \geq \omega$. Then $E$ has an $\omega$-good coloring over
$X$.
\et

\bp
The statement follows from Theorem \ref{asse}.\ref{small} with $\kappa = \omega$, $\mu = |X|$, and $\gotH = E \cap [X]^{\leq 2}$ counted with multiplicity.
\ep

From now on we work for $n$-good colorings with $n <\om$. The case $n= 1$ is trivial, so we start with
$n= 2$. We have the following characterization for the existence of $2$-good colorings.

\bt
\lab{t:graph2}
Let $G = (V,E)$ be a graph. Then $E$ has a $2$-good coloring if and only if no connected component of $G$ is an odd cycle.
\et
\bp
It is easy to see that the condition is necessary, as an odd cycle has no $2$-good coloring.

To prove sufficiency, observe that we can assume $G$ is connected
since it is sufficient to color the 
connected components separately.

\bl\label{UGL}
There exists a subgraph $G' = (V',E')$ of
$G$ and a coloring $c' \colon E'\to 2$ such that   $V' \neq \es$ and  for every $v\in V'$, $d_{G}(v)\ge 2$ implies $c'[E'(v)]=2$.
\el
\bp
If there is $x\in V$ with $d_G(x)=1$ then 
$V'=\{x\}$ and $E'=\empt$ works.
So we can assume that $d_G(x)\ge 2$ for each $x\in V$. Depending on the subgraphs of $G$ we
distinguish several cases.

\emph{Case I: $G$ contains an even cycle, or a path which is
infinite in both directions.} Since even cycles and paths
infinite in both directions have 2-good colorings we can choose $G'$ to be one
of these subgraphs. Note that a pair of multiple edges is an even cycle.

From now on we assume that $G$ contains no such subgraphs. Pick an arbitrary vertex $v
\in V$ and start a path from $v$ until it first fails to be vertex-disjoint.

\emph{Case II: we get an infinite path in this direction.} Let us start another path
from $v$ until it is vertex-disjoint from itself and from the previous infinite path. As we have no doubly infinite
paths, the second path has to terminate, say at $w \in V$. We obtained a cycle on $w$ and
an infinite
path starting from $w$. Let us call such a configuration an \emph{infinite
lasso}. It is easy to check that we get a 2-good coloring of our infinite
lasso if we color its edges by the alternating coloring in such a way that we start the coloring at $w$ and first we color the edges of the cycle on $w$. Thus in this case we can set $G'$ to be an infinite lasso.

\emph{Case III: our (first) path from $v$ reaches a vertex visited before.}
Since $G$ contains no even cycles we get an odd cycle $C$. Since $G$ is not an odd cycle, there is a vertex $w$ of $C$ with
$d_{G}(w) \geq 3$. Let us start a path from $w$ disjoint from $C$.
If we get an infinite lasso then we are done by Case II. Otherwise the path
reaches either a vertex of $C$ or a vertex of the path itself. If it reaches $C$
then it has to reach it at $w$: else $G$ would contain an even
cycle since for two odd cycles intersecting each other in a finite path, removing the intersection results an even cycle.

Hence we obtain two disjoint cycles connected by a path, possibly of
length $0$. Let $G'$ be this graph. As for the infinite lasso, color the edges of this graph the
alternating way, starting from $w$ and coloring first a cycle containing
$w$. It is easy to check that this is a 2-good coloring of $G'$.
\ep

Now we go back to the proof of Theorem \ref{t:graph2}.  For an ordinal $\xi$
to be specified later, we define a
sequence of partial colorings $c_{\alpha} \colon E \rar 2$ $(\alpha < \xi)$ such that
\ben[(i)]
\item\label{2i} $\dom(c_{\alpha}) \subsetneq \dom (c_{\alpha'})$ and $c_{\alpha'}|_{\dom
(c_{\alpha})} =c_{\alpha}$ $(\alpha < \alpha' <
\xi)$,
\item\label{2ii} if $v\in V[ \dom (c_{\alpha})]$ and $d_{G}(v)\geq 2$ then $ c_{\alpha}[E(v)] = 2$
$(\alpha < \xi)$,
\item\label{2iii} $V= V\left[\bigcup _{\alpha < \xi} \dom (c_{\alpha})\right]$.
\een

Once this done the function $c \colon E \rar 2$, $c = \bigcup _{\alpha < \xi}
c_{\alpha}$ is a 2-coloring of $E$ by (\ref{2i})  which is 2-good by (\ref{2ii}) and (\ref{2iii}).

To start the construction, by Lemma \ref{UGL} 
we have a partial coloring  $c'$ which works as $c_0$.
 Let $\alpha$ be an ordinal and suppose
that $c_{\beta}$ is defined for every $\beta < \alpha$. If $V= V\left[\bigcup _{\beta <
\alpha} \dom (c_{\alpha})\right] $ set $\xi = \alpha$ and the construction is
done. Else set $c^{-}_{\alpha} =\bigcup _{\beta <
\alpha} c_{\beta}$. We have $V \sm V[\dom (c^{-}_{\alpha})] \neq \es$. As $G$ is connected, there
exists an edge $\{u,v\} \in E \sm \dom (c^{-}_{\alpha})$ such that $u \in V[\dom (c^{-}_{\alpha})]$ and $v \in V \sm V[\dom (c^{-}_{\alpha})]$.
Start a path $P$ from $u$ whose first edge is $\{u,v\}$, and we keep extending $P$ as long as it
is  edge disjoint from $\dom (c^{-}_{\alpha})$.
This $P$ can be infinite, or it can end either in $V[\dom (c^{-}_{\alpha})]$ or in
a vertex of $P$ or in a vertex $x\in V$ with $d_G(x)=1$.
Let $c_{\alpha}$ be the partial coloring extending $c^{-}_{\alpha}$ where we
color the edges of $P$ the alternating way starting with the edge $\{u,v\}$. Then
$c_{\alpha}$ clearly satisfies (\ref{2i}); we have to show that $c_{\alpha}$ also satisfies  (\ref{2ii}). To this end, let
$w\in V[ \dom (c_{\alpha})]$ with $d_{G}(w)\geq 2$. By the definition of $P$, either $w\in V[\dom(c^{-}_{\alpha})]$ or $d_{P}(w) \geq 2$.
If $w\in V[\dom(c^{-}_{\alpha})]$ then $ c^{-}_{\alpha}[E(w)] = 2$
by the inductive assumption (\ref{2ii}); while if $d_{P}(w) \geq 2$ then 
$c_{\alpha}[E(w)]=2$ because we color the edges of $P$
the alternating way. 

Since $\dom(c_{\alpha})$ $(\alpha < \xi)$ are strictly increasing, this
transfinite procedure terminates at some ordinal $\xi$. The resulting sequence
$(c_{\alpha})_{\alpha < \xi}$ satisfies (\ref{2i})-(\ref{2iii}), so the proof is complete.
\ep

Clearly, an $n$-regular graph has an $n$-good coloring if and only if  its edge
chromatic number is $n$. It is a well-known theorem of Vizing that the
edge chromatic number of a simple finite graph is either $\Delta(G)$
or $\Delta(G)+1$ (see e.g.\ \cite[Theorem 28.2 p.\ 467]{Sc}). But  to decide e.g.\ whether a 3-regular graph is 3-chromatic or not
is an NP-complete problem (see e.g.\ \cite[Theorem 28.3 p.\ 468]{Sc}).
Hence we cannot hope for a very simple analogue of Theorem
\ref{t:graph2} for $n \geq 3$.

It remains to extend the Theorem of R.\ P.\ Gupta to infinite graphs.
\bt\label{Gupta_inf}
Let $1 \le n <\om$.
Let $G = (V,E)$ be a graph and let $X \su V$ be such that for every $x \in X$ we have $d_{G}(x)  \geq n+\mu(G)$. Then $E$ has an $n$-good coloring over $X$.
\et

\bp
For finite graphs this is Theorem \ref{Gupta} due to Gupta.
If $G$ is locally finite, i.e.\ all degrees are finite, then an easy
compactness argument yields the result. If $\mu(G) \geq \omega$ we are done by Theorem \ref{UJU}.

If $G$ is arbitrary with $\mu(G) < \omega$, we
construct a locally finite graph $G' = (V',E')$, as follows. Set $V_{0} = \{v \in V \colon d_{G}(v) < \omega\} $ and $V_{1} = \{v \in V \colon d_{G}(v) \geq \omega\} $. Let $V' = V_{0} \cup \{\langle v,\alpha\rangle \colon v \in V_{1},~\alpha < d_{G}(v)\}$. 
We define $E'$ by setting $G'[V_{0}] = G[V_{0}]$ and by distributing for $v \in V_{1}$, using $d_{G}(v) = d_{G}(v) \times (n+\mu(G))$, the edges $E(v)$ onto $(\langle v,\alpha\rangle)_{\alpha < d_{G}(v)}$ ``uniformly": that is, $E'$ is constructed in such a way that
\ben[(i)]
\item for every $v\in V_{0}$, $d_{G'}(v) = d_{G}(v)$,
\item for every $v \in V_{1}$ and $\alpha < d_{G}(v)$, $d_{G'}(\langle v,\alpha\rangle) = n+\mu(G)$,
\item for every $v \in V_{1}$ and $\alpha < d_{G}(v)$, $E'(\langle v,\alpha\rangle)$ contains no multiple edges. 
\een
This is clearly possible, and we have $\mu(G') \leq \mu(G)$. Set $X' = (X \cap V_{0}) \cup \{\langle v,0\rangle \colon v \in X \cap V_{1}\}$. We have that  $G'$ is locally finite and $d_{G'}(x) \geq n+\mu(G)$ hence $d_{G'}(x) \geq n+\mu(G')$  $(x \in X')$. So there is an $n$-good coloring of $G'$ over $X'$. By merging, for every $v \in V_{1}$, the vertices $\langle v,\alpha\rangle \in V'$ $(\alpha < d_{G}(v))$ to one vertex we get a graph isomorphic to $G$ and $X'$ is mapped onto $X$.
Thus the $n$-good coloring of $G'$ over $X'$ yields an $n$-good coloring of $G$ over $X$.
\ep

\section{Intervals in linearly ordered sets}
\lab{s:in}

  Let $\mc{L} = (L, \leq)$ be a linearly ordered set and let $\conv(\mc{L})$ denote the family of convex subsets of $L$.
In this section we prove the following two results. The first establishes maximally good coloring for convex covers.

\bt\label{tm:inf}
Let $(L,\le) $ be an ordered set and let $\gotH$ be a cover of $L$ with $\iH \subseteq \conv(L)$. Then $\gotH$ has a maximally good coloring.
\et

The second gives the splitting of $k$-covers.

\bt\label{tm:fin}
Let $(L,\le) $ be an ordered set and, for some (finite or infinite) cardinal $k$, let $\gotH$ be a $k$-fold cover of $L$ with $\iH \subseteq \conv(L)$. Then $\gotH$ has a $k$-good coloring.
\et

As we noted in the introduction, Theorem \ref{tm:fin} was obtained in \cite{AHM} much earlier than our investigations. In spite of \cite{AHM}, we decided to treat the splitting of convex covers of linearly ordered sets in the present paper because Theorem \ref{tm:inf} is new, and we found a significantly simpler proof of Theorem \ref{tm:fin} than the one in \cite{AHM}.

The heart of the proof of Theorem \ref{tm:inf} is the following general statement on maximally good colorings.

\bt
\label{tm:ex}
Let $X$ be a set and let $\gotH$ be a simple cover on $X$.
If for each $\gotK \subseteq \gotH$ there is $\gotJ \subseteq \gotK$
such that
\begin{enumerate}[(e1)]
\item $\cup\gotJ=\cup\gotK$,
\item $\gotJ$ has a maximally good coloring,
\end{enumerate}
then $\gotH$ has a maximally good coloring.
\et
\bp Let $\{ H_{\alpha}\colon \alpha < |\gotH|\}$ be an enumeration of $\gotH$.
By transfinite recursion on ${\alpha}<|\gotH|$ we define families
$\gotJ_{\alpha}\subseteq \gotH$ satisfying $H_{\alpha}\in \bigcup_{\beta \leq \al} \gotJ_{\beta}$, and maximally good colorings
$c_{\alpha}:\gotJ_{\alpha}\to \gotO$, as follows.
Let ${\alpha}<|\gotH|$ be arbitrary and suppose $\gotJ_{\nu}$ and $c_{\nu}$ are constructed for
${\nu}<{\alpha}$. Set
$\gotK_{\alpha}=\gotH\setm \cup\{\gotJ_{\nu}:{\nu}<{\alpha}\}$.
Since $\gotK_{\alpha}\subs \gotH$, by assumption
we have a family $\gotJ_{\alpha}\subs\gotK_{\alpha}$
with $\cup\gotJ_{\alpha}=\cup\gotK_{\alpha}$
and  a maximally good coloring $c_{\alpha}$
of $\gotJ_{\alpha}$.
 If $H_{\alpha} \notin \bigcup_{\beta \leq \al} \gotJ_{\beta}$, we put $H_{\alpha}$ into $\gotJ_{\alpha}$ and we set $c_{\alpha}(H_{\alpha}) = 0$. So we can assume $H_{\alpha} \in \bigcup_{\beta \leq \al} \gotJ_{\beta}$. This completes the $\al^{\textrm{th}}$ step of the construction. Then we have
$\gotH=\bigcup_{{\alpha}<|\gotH|}\gotJ_{\alpha}$.

Let $c:\gotH \to \gotO $,  $c( H)={\alpha}+c_{\alpha}(H)$ for $ H\in\gotJ_{\alpha}$ $(\alpha < |\gotH|)$.
By $\gotK_{0} =\gotH$ we have $\cup \gotJ_0 = \cup\gotH$ and $c|_{ \gotJ_0}= c_0$, so $c$ is 1-good.

Before proving that $c$ is $[\omega,\infty)$-good, let us observe that if $x \in X$ and $\alpha$ is an ordinal such that  $\gotJ_{\beta}(x) \neq \es$ for each ${\beta}<{\alpha}$, then
$0\in c_{\beta}[\gotJ_{\beta}(x)]$ and so ${\beta}\in c[\gotJ_{\beta}(x)]$ $({\beta}<{\alpha})$.
Hence ${\alpha}\subs c[\gotH(x)]$.

To see that $c$ is $[\omega,\infty)$-good, pick $x \in L$ and suppose ${\kappa}=|\gotH(x)|\ge {\omega}$. We distinguish several cases. If $\gotJ_{\beta}(x) \neq \es$ for each ${\beta}<{\kappa}$, then by the previous observation ${\kappa}\subs c[\gotH(x)]$, as required.

So suppose $\gotJ_{\beta}(x) = \es $ for some
${\beta}< {\kappa}$; fix a minimal such $\beta$.
Then $\gotH(x)=\bigcup\{\gotJ_{\alpha}(x):{\alpha}<{\beta}\}$.
Thus for each cardinal ${\lambda}<{\kappa}$ there is an
${\alpha(\lambda)}<{\beta}$ such that $|\gotJ_{\alpha(\lambda)}(x)|\ge \max\{\omega, \lambda^{+}\}$.
Then $\gotJ_{\gamma}(x)\ne \empt $ $({\gamma}<{\alpha(\lambda)})$ and so
$ {\alpha(\lambda)}\ss c[\gotH(x)]$ by our observation. Moreover,
$\max\{\omega, \lambda^{+}\} \ss c_{\alpha(\lambda)}[\gotJ_{\alpha(\lambda)}]
$ and so $[\alpha(\lambda), {\alpha(\lambda)}+{\lambda}^+ )\ss c[\gotH(x)]$, as well. By putting these together we obtain ${\alpha(\lambda)}+{\lambda}^+ \ss c[\gotH(x)]$, so since ${\kappa}=\sup \{{\lambda}^+:{\lambda}<{\kappa}\}$
we concluded ${\kappa} \ss c[\gotH(x)]$, as required.
\ep

In the following two lemmas, for any linearly ordered set $\mc{L}$, we establish the existence of a maximally good coloring for special subfamilies of $\conv(L)$.
 For $x\in L$ set
$(-\infty, x]=\{y\in L:y\le x\}$ and $[x,+\infty)=\{y\in L: x\le y\}$.
We define $$\tail(\lcal)=\{I\subs L:[x,+\infty)\subs I\text{ for
    each $x\in I$}\}.$$ Clearly, $\tail(\lcal)\subsetneq\conv(\lcal)$ provided $|L|\geq2$.

\bl\label{lm:ord}
Every simple cover $\gotH$ with $\iH \subs \tail(\lcal)$
has a maximally good coloring.
\el

\bp
We intend to apply Theorem \ref{tm:ex}. To this end, it is enough to show that for every $\gotK \ss \tail(\lcal)$ there is a $\gotJ \ss \gotK$ satisfying $\cup\gotK = \cup \gotJ$ such that $\gotJ$ has a maximally good coloring.

Let $\gotK \subs \tail(\lcal)$ be a cover. For some regular cardinal ${\kappa}$
there is a strictly increasing chain
$\gotJ=\{J_{\nu}:{\nu}<{\kappa}\}$ of elements of $\gotK$
such that $\cup\gotJ=\cup\gotK$. Note that we may have $\kappa = 1$.

Let $f:{\kappa}\to {\kappa}$ be a ${\kappa}$-abundant map, i.e.\ for every $\lambda < \kappa$ we have $|f^{-1}(\lambda)| = \kappa$. Define $c: \gotJ\to \gotO$ by $c(J_{\lambda})=f(\lambda)$ $(\lambda < \kappa)$. Clearly, $c$ is a maximally good coloring of $\gotJ$, which completes the proof.\ep

\bl\label{lm:point}
Fix $a\in L$. Let $\gotH$ be a simple cover with $ \iH \subs \conv(\lcal)(a)$. Then $\gotH$ has a maximally good coloring.
\el

\bp
Again, we intend to apply Theorem \ref{tm:ex}, thus it is enough to show that for every $\gotK \ss \conv(\lcal)(a)$ there is a $\gotJ \ss \gotK$ satisfying $\cup\gotK = \cup \gotJ$ such that $\gotJ$ has a maximally good coloring.

So let  $\gotK\subs \conv(\lcal)(a)$ be a cover. By the definition of maximally good coloring we can assume $\cup \gotK=L$. For some regular cardinal ${\kappa}$
there is a  family
$\gotJ^+=\{J_{\nu}:{\nu}<{\kappa}\}$ of elements of $\gotK$
such that
$\{J_{\nu}\cap [a,+\infty):{\nu}<{\kappa}\}$ is strictly increasing
and $\gotJ^+$ covers $[a,+\infty)$. Note that we may have $\kappa=1$. We can apply  Lemma \ref{lm:ord} for $\gotJ^+$ as a cover over $(-\infty,a]$ to obtain a maximally good coloring $c:\gotJ^+\to \gotO$  of $\gotJ^+$ over $(-\infty,a]$. 

Let $f:{\kappa}\to {\kappa}$ be a ${\kappa}$-abundant map.
Define $h:{\kappa}\to {\kappa}$ by
\begin{displaymath}
h({\beta})=\left\{
\begin{array}{ll}
f({\xi}+n)&\text{if ${\beta}={\xi}+2n+1$ for some ${\xi}\in \Lim({\kappa})$},\\
{\xi}+n&\text{if ${\beta}={\xi}+2n$ for some ${\xi}\in \Lim({\kappa})$},
\end{array}
\right .
\end{displaymath}
and let $d^+:\gotJ^+\to \gotO$, $d^+(J_{\nu})=h(c(J_{\nu}))$ $(\nu < \kappa)$.
Then $d^+$ is a maximally good coloring of $\gotJ^+$.

If $\gotJ^+$ covers $L$, $\gotJ=\gotJ^+$ satisfies the
requirements. If not, take $\gotK'=\gotK\setm \gotJ^{+}$.
Then $\gotK'$ covers $(-\infty,a]$, so by repeating the previous argument for
$(-\infty,a]$ instead of $[a,+\infty)$
we can find a family $\gotJ^-\subs \gotK'$ covering
$(-\infty,a]$ with a maximally good coloring
$d^-$.

Put $\gotJ=\gotJ^+\cup\gotJ^-$ and $d =d^-\cup d^+$.
Then $\cup\gotJ=L$ and $d$ is a maximally good coloring of $\gotJ$, which completes the proof.\ep

\medskip

\textbf{Proof of Theorem \ref{tm:inf}.} By Proposition \ref{mred} we can assume $\gotH$ is simple.
By Theorem \ref{tm:ex} it is enough to prove that for every cover $\gotK$ with $\gotK \ss \conv(\mc{L})$
there is a subfamily $\gotJ\subs \gotK$
such that
\begin{enumerate}[(a)]
\item $\cup\gotJ=\cup\gotK$,
\item $\gotJ$ has a maximally good coloring $c_{\gotJ}$.
\end{enumerate}

Consider the equivalence relation $R$ on $L$ generated by the relation
$\bigcup\{I\times I:I\in \gotK\}$. The equivalence classes of $R$ give a partition of $L$ and every $I\in \gotK$ is contained in
some equivalence class. Hence we can construct $\gotJ$ and $c_{\gotJ}$ for each
equivalence class separately. Therefore we can assume
we have only one equivalence class. Hence for every $z^{-}, z^{+} \in L$, $[z^{-}, z^{+}]$ can be covered by  finitely many members of $\gotK$.

 Let $z\in L$ be arbitrary.

\bc\label{korok}
If
\begin{itemize}
\item[$(\circ)$] for each $x\in [z, +\infty)$ there is $ y\in [z, +\infty)$
 such that $\cup\gotK(x)\subs (-\infty,y]$
\end{itemize}
then
there is $\gotJ^+\in \br \gotK;{\omega};$ such that
\begin{itemize}
\item[$(\circ\circ)$] $\gotJ^+$ covers
$[z,+\infty)$ and $|\gotJ^+(x)|<{\omega}$
for each $x\in L$
\end{itemize}
\ec
\bp
We define recursively a partition
$\{L(n)\colon n\in {\omega}\}$  of $[z,+\infty)$ by setting
$L(0)=\{z\}$, and for $0 < n < \omega$,
\begin{displaymath}
L(n)=\{y\in [z,+\infty)\colon I\cap L(n-1)\ne \empt
\text{ for some $I\in \gotK(y)$} \}\setm \textstyle
\bigcup_{k < n} L(k).
\end{displaymath} Since $L$ is one equivalence class of $R$, $[z, +\infty) = \bigcup_{n < \omega}L(n)$, indeed. Note that some $L(n)$ can be empty, e.g.\ if $L$ has a maximal element.

We show that for each $n < {\omega}$ there is  an $\gotI_n\in \br \gotK;\le 2;$ such
 that $\gotI_n$ covers $L(n)$ and $I \cap L(n) \neq \es$ $(I \in \gotI_{n})$. This is obvious if $n=0$ or $L(n) = \es$. If $n \neq 0$, $L(n) \neq \es$ but $L(n+1) = \es$ then by $(\circ)$, $L$ has a maximal element $m$ and $m \in L(n)$. By definition, there is an $I \in \gotK(m)$ with $I \cap L(n-1) \neq \es$, so $\gotI_{n} = \{I\}$ fulfills the requirements. Finally if $n \neq 0$, $L(n) \neq \es$ and $L(n+1) \neq \es$ then pick a $y \in L(n+1)$. By definition, there is an $I \in \gotK$ with $y \in I$ and $I \cap L(n) \neq \es$. Let $y' \in I \cap L(n)$ and let $I' \in  \gotK$ with $y' \in I'$ and $I' \cap L(n-1) \neq \es$. Then $\gotI_{n} = \{I, I'\}$ fulfills the requirements.

 Let $\gotJ^+=\cup\{\gotI_n \colon n<{\omega}\}$. Since $[z, +\infty) = \bigcup_{n < \omega}L(n)$, $\gotJ^+$ covers $[z,+\infty)$. Observe that
 for each $n < \omega$, if $x \in L(n)$ then $I\in \gotI_{n+2}$ implies $x \notin I$. Hence  $|\gotJ^+(x)|<{\omega}$ $(x \in L)$, as required.
\ep

Let us return to the proof of Theorem \ref{tm:inf}.
If $(\circ)$ holds then let $\zp=z$ and fix a family
$\gotJ^+\in \br \ical;{\omega};$ satisfying $(\circ\circ)$.
Otherwise pick $\zp\in [z,+\infty)$ such that $\gotK(\zp)$ covers $  [\zp,+\infty)$
and let $\gotJ^+=\gotK(\zp)$.

By applying Proposition \ref{korok} to $L$ with reversed order, we can show that if
\begin{itemize}
\item[$(\diamond)$] for each $x\in (-\infty,z]$ there is $ y\in (-\infty,z]$
 such that $\cup\gotK(x) \subs [y,+\infty)$
 \end{itemize}
 then there is
 a family
 $\gotJ^-\in \br \gotK;{\omega};$ such that
\begin{itemize}
\item[$(\diamond\diamond)$] $\gotJ^-$ covers
$(-\infty,z]$ and $|\gotJ^-(x)|<{\omega}$
for each $x\in L$.
\end{itemize}
If $(\diamond)$ holds let $\zm=z$ and fix a family  $\gotJ^-$ satisfying $(\diamond\diamond)$.
Otherwise pick $\zm\in (-\infty,z]$ such that $\gotK(\zm)$ covers
  $(-\infty,\zm]$
and let $\gotJ^-=\gotK(\zm)$.
Finally pick  $\gotJ^{0}\in \br \gotK;<{\omega};$ which covers
$[\zm,\zp]$. Let $\gotJ=\gotJ^{-} \cup \gotJ^{0}\cup \gotJ^+$.
Then $\gotJ$ covers $L$.

The families $\gotJ^+$, $\gotJ^-\setm \gotJ^+$, and $\gotJ^{0}\setm
(\gotJ^+\cup \gotJ^-)$ have maximally good colorings $c^+$, $c^-$ and $c^0$ respectively,
because they are either ``locally finite'' or Lemma \ref{lm:point}
can be applied.  Thus $c_{\gotJ}=c^+\cup c^-\cup c^0$ is a maximally good coloring of $\gotJ$.
\ep

We close this section with the proof of Theorem \ref{tm:fin}.

\medskip

\textbf{Proof of Theorem \ref{tm:fin}.}
If $k$ is an infinite cardinal the statement follows immediately from Theorem \ref{tm:inf}. So let $k < \omega$; we prove the statement
by induction on $k$. For $k=1$ the statement is trivial.

Let $k \geq 2$ and suppose the theorem is true for $k-1$. As in the proof of Theorem \ref{tm:inf}, consider the equivalence relation $R$ on $L$ generated by the relation
$\bigcup\{H\times H \colon H\in \hcal\}$. The equivalence classes of $R$ give a partition of $L$ and every $H\in \hcal $ is contained in
some equivalence class, hence we can construct the $k$-good coloring of $\gotH$ for each
equivalence class separately. Therefore we can assume that
we have only one equivalence class.

\bc\label{split}
Let $I\subs L$ be a convex set and $y\in I$.
If $\gotH$ has a $k$-good coloring over  $I\cap (-\infty,y]$
and another $k$-good coloring over  $I\cap [y,+\infty)$ then it has a $k$-good coloring over  $I$, as well.
\ec
\bp Fix two $k$-good colorings $c_{-} \colon \gotH \rar k$ and $c_{+} \colon \gotH \rar k$ over $I\cap (-\infty,y]$ and $I\cap [y,+\infty)$, respectively. By thinning out the domain of $c_{-}$ we can assume that for each $i<k$ the family $[c_{-}^{-1}(i)](y)$ has an enumeration
$\{J^i_{-}(\gamma)\colon \gamma<{\kappa}_i\}$ for some regular cardinal
${\kappa}_i$ such that
$\{J^i_{-}(\gamma)\cap (-\infty,y]:{\gamma}<{\kappa}_i\}$ is strictly
increasing and so
for each cofinal subset $\Gamma\subs {\kappa}_i$
the family $(c_{-}^{-1}(i) \setm [c_{-}^{-1}(i)](y))\cup\{J^i_{-}(\gamma) \colon{\gamma}\in \Gamma\}$
covers $I\cap (-\infty,y]$.
Let us remark that ${\kappa_{i}}$ can be finite, namely $1$.

 Similarly,
we can thin out the domain of $c_{+}$ such that for each $i<k$
the family $[c_{+}^{-1}(i)](y)$ has an enumeration
$\{J^i_{+}(\gamma)\colon \gamma<{\lambda}_i\}$ for some regular cardinal
${\lambda}_i$ such that for each cofinal subset $\Gamma\subs {\lambda}_i$
the family $(c_{+}^{-1}(i) \setm [c_{+}^{-1}(i)](y))\cup\{J^i_{+}(\gamma) \colon {\gamma}\in \Gamma\}\}$
covers $I\cap [y,+\infty)$.

Then by passing to cofinal subsets of $[c_{-}^{-1}(i)](y)$ and $[c_{+}^{-1}(i)](y)$ we can assume that for each $i,j<k$
if $[c_{-}^{-1}(i)](y)\cap [c_{+}^{-1}(j)](y) \ne \empt$ then
${\kappa}_i={\lambda}_j=1$ and so $[c_{-}^{-1}(i)](y)=[c_{+}^{-1}(j)](y)$. So there is a bijection $f:k\to k$ such that if $[c_{-}^{-1}(i)](y)\cap [c_{+}^{-1}(j)](y)\ne \empt$ then $j=f(i)$.

Define $c \colon \gotH \rar k$ by $c(\langle H, \alpha \rangle) = i$ if  $c_{-}(\langle H, \alpha \rangle) =i$ or $c_{+}(\langle H, \alpha \rangle) = f(i)$ $(\langle H, \alpha \rangle \in \gotH)$.  The definition of $c$ is valid and $c$ is a $k$-good coloring of $\gotH$ over $I$. This completes the proof.
\ep

Define the relation $\equiv$ on $L$ by
$x\equiv y$ if and only if there exists a $k$-good coloring of $\gotH$ over $[x,y]$. By Proposition \ref{split}, $\equiv$ is an equivalence relation on $L$.  Moreover, we have the following.

\bc\label{5.8} For every $H\in \hcal$, $H$ is contained in one
equivalence class of $\equiv$.
\ec
\bp Let $H\in\hcal$ and $\{x,y\}\in \br H;2;$.
Then $\gotH\setm \{\langle H,0\rangle\}$ is a $(k-1)$-fold cover of
$L$. Hence by the inductive hypothesis,
$\gotH\setm \{\langle H,0\rangle\}$ has a $(k-1)$-good coloring  $c \colon \gotH\sm \{\langle H,0\rangle\} \rar k-1$ over $[x,y]$.
Extend $c$ by setting $c(\langle H,0\rangle)= k$; then $c$ is a
$k$-good coloring over $[x,y]$.
\ep

\bc\label{E}
Let $E$ be an equivalence class of $\equiv$. Then there is a $k$-good coloring of $\gotH$ over $E$.
\ec
\bp
Take an arbitrary $y\in E$. Since $E$ is convex, by Proposition \ref{split} it is enough to prove that $\gotH$ has a $k$-good coloring over $E\cap [y,+\infty)$ and over $E\cap (-\infty, y]$. We prove only that $\gotH$ has a $k$-good coloring over $E\cap [y,+\infty)$, the proof of the other statement is similar. We distinguish several cases.

Suppose first that there is
$H\in \iH$ such that $H$ is cofinal in $E$. Fix $z\in H\cap E \cap [y,+\infty)$; then $[z, +\infty) \ss [y,+\infty)$.
So $\gotH\setm \{\langle H,0\rangle\}$ is a $k-1$-fold cover of $E \cap [z,+\infty)$. So by the inductive hypothesis,
there is $c \colon \gotH \sm \{\langle H,0\rangle\} \rar k-1$, a $k-1$-good coloring of $\gotH \sm \{\langle H,0\rangle\} $ over $E\cap[z,+\infty)$. Then extending $c$ by setting $c(\langle H,0\rangle)=k$ yields a $k$-good coloring of $\gotH$ over $E \cap [z,+\infty)$.
Since $y \equiv z$, $\gotH$ has a $k$-good coloring over $[y,z]$, so by Proposition \ref{split} we have that $\gotH$ has a $k$-good coloring over $E\cap [y,+\infty)$, as well.

From now on assume that there is no  $H\in \iH$ such that  $H$ is cofinal in $E\cap [y,+\infty)$. If there is $z\in [y,+\infty)$ such that
$\cup\iH(z)$ is cofinal in $E$, then since for $H \in \iH$, $H$ is not cofinal in $E$, $\gotH$ has a $k$-good coloring over
$E\cap [z,+\infty)$.
Since  $y \equiv z$, $\gotH$ has a $k$-good coloring over $[y,z]$.
So by Proposition \ref{split}, $\gotH$ has a $k$-good coloring over
$E\cap [y,+\infty)$, as well.

In the sequel we assume in addition that for every $z\in E\cap [y,+\infty)$, $\cup\iH(z)$ has an upper bound in $E$. We define recursively a strictly increasing sequence
$(x_{n})_{n < \omega}\subs E\cap [y,+\infty)$, as follows. Let $x_{0} = y$. If $0 < n < \omega$ and $x_{n-1}$ is already defined, let $b_{n-1}$ be an upper bound of $\cup\gotH(x_{n-1})$, and let $x_{n}$ be an upper bound of $\cup\gotH(b_{n-1})$.  Then $\cup \gotH(x_{n-1}) \ss (-\infty,b_{n-1}]$ and $\cup\gotH(x_{n}) \ss (b_{n-1},+\infty)$ imply $$\gotH(x_{n})\cap \gotH(x_{n'})=\empt~(  n<n'<{\omega}),$$
and by our assumption that $L$ is one equivalence class of $R$,
 $\{x_{n}\colon n<{\omega}\}$ is cofinal in $E$.

For every $n<{\omega}$ we have $x_{n} \equiv x_{n+1}$ so there is $c_{n} \colon \gotH \rar k$, a $k$-good coloring of $\gotH$ over   $[x_n,x_{n+1}]$. Fix $n < \omega$; by thinning out the domain of $c_{n}$ we can assume that for each $i<k$ the family $[c_{n}^{-1}(i)](x_{n})$ has an enumeration
$\{J^i_{n}(\gamma)\colon \gamma<{\kappa}^{i}_n\}$ for some regular cardinal
${\kappa}^{i}_n$ such that
$\{J^i_{n}(\gamma)\cap [x_{n},x_{n+1}] \colon{\gamma}<{\kappa}^{i}_{n}\}$ is strictly
increasing, and so
for each cofinal subset $\Gamma\subs {\kappa}^{i}_{n}$
the family $(c_{n}^{-1}(i) \setm [c_{n}^{-1}(i)](x_{n}))\cup\{J^i_{n}(\gamma):{\gamma}\in \Gamma\}$
covers $[x_{n},x_{n+1}]$. Similarly, we can assume that for each $i<k$ the family $[c_{n}^{-1}(i)](x_{n+1})$ has an enumeration
$\{B^i_{n}(\gamma)\colon \gamma<{\lambda}^{i}_n\}$ for some regular cardinal
${\lambda}^{i}_n$ such that
$\{B^i_{n}(\gamma)\cap [x_{n},x_{n+1}] \colon{\gamma}<{\lambda}^{i}_{n}\}$ is strictly
increasing, therefore
for each cofinal subset $\Gamma\subs {\lambda}^{i}_{n}$
the family $(c_{n}^{-1}(i) \setm [c_{n}^{-1}(i)](x_{n+1}))\cup\{B^i_{n}(\gamma):{\gamma}\in \Gamma\}$
covers $[x_{n},x_{n+1}]$.

Then for every $n < \omega$ and $i < k$, we can pass to cofinal subsets of $[c_{n}^{-1}(i)](x_{n+1})$ and $[c_{n+1}^{-1}(i)](x_{n+1})$ in such a way that for each $i,j<k$
if $[c_{n}^{-1}(i)](x_{n+1}) \cap [c_{n+1}^{-1}(j)](x_{n+1}) \ne \empt$ then
${\lambda}^{i}_{n}={\kappa}^{j}_{n+1}=1$ and so $[c_{n}^{-1}(i)](x_{n+1})=[c_{n+1}^{-1}(j)](x_{n+1})$. So there is a bijection $f_{n}:k\to k$ such that if $[c_{n}^{-1}(i)](x_{n+1})\cap [c_{n+1}^{-1}(j)](x_{n+1})\ne \empt$ then $j=f_{n}(i)$. Write $g_{0} = \textrm{Id}$ and $g_n=f_{n-1}\circ f_{n-2}\circ \dots \circ f_0$ $(0 < n < \omega)$.

For every $\langle H, \alpha \rangle \in \gotH$ and $i<k$ define $c(\langle H, \alpha \rangle) = i$ if and only if for some $n  < \omega$, $\langle H, \alpha \rangle \in \dom(c_{n})$ and $c_{n}(\langle H, \alpha \rangle) = g_n(i)$. This definition makes sense and $c$ is a $k$-good coloring of $\gotH$ over $E \cap [y,+\infty)$, which completes the proof.
\ep

We are ready to complete the proof of Theorem \ref{tm:fin}. By assumption, $L$ is one equivalence class of the relation $R$. So by Proposition \ref{5.8}, $L$ is one equivalence class of $\equiv$. Therefore by Proposition \ref{E} there is a $k$-good coloring of $\gotH$, which finishes the proof. \ep

\section{Closed sets}\label{closedsets}

Towards the investigation of splitting of covers with special geometric properties let us tackle closed covers, i.e.\ that variant of the problem where the sets in the cover are closed. The study of this special case is motivated by the facts that, apart from considering open covers, this is the simplest topological constraint one can impose; even for closed covers we get independence of ZFC by very strong means; these results will be very useful for treating the problem of covers by compact convex sets.

Obviously, we have to specify the topological spaces where closed covers are considered. Observe that similarly to the proof of Theorem \ref{MA} below, the construction of Theorem \ref{asse}.\ref{power} can be carried out in such a way that the covering sets $H_{\alpha}$ are closed in $2^{\kappa}$ endowed with the product topology. Since our purpose is not to find suitable topologies for general constructions but to establish independence of ZFC for natural topological spaces, in this section we restrict our attention to covers of $\real$, or equivalently to covers of $\omega^{\omega}$ and $2^{\omega}$.

As we shall see in Proposition \ref{Marci}, if $\gotH$ is a closed cover of $\real$ and $|\gotH| < \cov(\mc{M})$ then $\gotH$ has a countable subcover. In particular, for $\omega < \kappa < \cov(\mc{M})$, a $\kappa$-fold closed cover \emph{of cardinality $\kappa$} has a $\kappa$-good coloring. There are models of ZFC where even Borel covers of special cardinalities of the real line satisfy a similar Lindel\"of like property. In \cite{M}, A.\ Miller showed that in a model obtained from a model of CH by adding $\omega_{3}$ many Cohen reals, every cover  of $\real$ by $\omega_{2}$ many Borel sets has an $\omega_{1}$ subcover. Here the corresponding splitting result says that if $\gotH$ is an $\omega_{2}$-fold Borel cover of $\real$ and $|\mc{H}| = \omega_{2}$ then $\gotH$ has an $\omega_{2}$-good coloring. However, these are very special settings as far as splitting is concerned, so we do not pursue our investigations in this direction. For more background on covering numbers related to closed sets see \cite{MS}.

In this section our main results are the following. In Theorem \ref{MA} we obtain that if $MA_{\kappa}(\sigma$-centered$)$ holds there exists a $\kappa$-fold closed cover of $\real$, consisting of translates of \emph{one} compact set, which cannot be partitioned into two subcovers. In particular, we obtain in ZFC that there exists an $\omega$-fold closed cover of $\real$, consisting of translates of one compact set, which cannot be partitioned into two subcovers. Finally in Theorem \ref{Cs_Cohen} we establish that in the Cohen real model every  closed cover of $\real$ has an $[\omega_{1}, \infty)$-good coloring. In this section $X$ denotes any of $\real$, $\omega^{\omega}$ or $2^{\omega}$; and $2^{\omega}$ is identified with $\mc{P}(\omega)$ the usual way.

\subsection{Martin's Axiom}\label{sMA}

This section is devoted to the following theorem.

\bt\label{MA} Let $\kappa$ be a cardinal satisfying $\omega \leq \kappa < 2^{\omega}$ and assume $MA_{\kappa}(\sigma\textrm{-centered})$. Then there exists a $\kappa$-fold simple closed cover of $X$ which cannot be decomposed into two disjoint subcovers. Moreover, in $\real$ the cover may consist of translates of one compact set.
\et

Since $MA_{\omega}(\sigma\textrm{-centered})$ holds in ZFC we obtain the following corollary.

\bcor\label{Cs_omega} There exists an $\omega$-fold closed cover of $X$ which cannot be partitioned into two subcovers. If $X = \real$ the cover can consist of translates of one compact set.
\ecor

\medskip

We prove Theorem \ref{MA} first in $\real$ since there we need to construct the cover using translates of one compact set. We fix some notation in advance. For a set $F \ss \real$ let $\<F\>_{\rac}$ denote the linear span of $F$ in $\real$ considered as a vector space over the rationals $\rac$. We set $\Sigma = 4^{\omega}$.

In order to construct a cover of $\real$ using translates of one compact set, we need the following auxiliary construction.

\bl\label{sov} 
There exist perfect sets $F,W\subs [0,1]$ and a sequence 
$(v_{n})_{n <  \omega} \ss [0,1]$ with 
$\lim _{n \rar \infty} v_{n} = 0$
such that  
\ben[(i)]
\item\label{ind000} $ F = \b\{ \sum_{i <
  \omega} \sigma(i)/4^{k_{i}} \colon \sigma \in \Sigma \j\}$
for some sequence $(k_{i})_{i < \omega} \ss \omega \sm \{0\}$;
\item\label{ind00} $W \cup \{v_{n} \colon n <
\omega\}$ is linearly independent over $\rac$;
\item \label{ind0}$\<F\>_{\rac} \cap \<W \cup \{v_{n} \colon n <
\omega\}\>_{\rac} = \{0\}$. 
\een
\el
\prf
 Let $(j_{i})_{i < \omega} \ss \omega\sm \{0\}$ satisfy
$j_{i+1}-j_i>i$ $(i < \omega)$. Write $k_i=j_{2i}$
and $\ell_i=j_{2i+1}$ $(i < \omega)$, and set $ F = \b\{ \sum_{i <
  \omega} \sigma(i)/4^{k_{i}} \colon \sigma \in \Sigma \j\}$
and $ U = \b\{ \sum_{i <
  \omega} \sigma(i)/4^{\ell_{i}} \colon \sigma \in \Sigma \j\}$. Then (i) holds.
  
  By \cite[Theorem 1 p.\ 141]{My},
there is a nonempty perfect set $U ' \subs U$
such that $U' \cup\{1\}$ is linearly independent over $\rac$; in particular, $\<U'\>_{\rac} \cap \rac = \{0\}$. Let $(w_n)_{ n<{\omega}}\subs U '$ be a strictly decreasing sequence. Set $v_n=w_n/(n+1)$ $(n < \omega)$ and let $W \ss U' \sm \{w_n \colon n<{\omega}\}$ be a nonempty perfect set. Then  
$\lim _{n \rar \infty} v_{n} = 0$, and (ii) holds.
  
It remains to verify (iii). First we show that $\<F\>_{\rac} \cap \<U\>_{\rac} \subs \rac$. To see this, for some $m,n < \omega$, let $f_{a} \in F$, $p_{a} \in \rac\sm \{0\}$ $(a < m)$ and $u_{b} \in U$, $q_{b} \in \rac\sm \{0\}$ $(b < n)$ satisfy \Keq\label{PO}\sum_{a < m} p_{a}f_{a} = \sum_{b < n} q_{b}u_{b}.\Zeq For every $a < m$ and $b < n$, let $\sigma_{a}, \tau_{b} \in \Sigma$ be such that $f_{a} = \sum_{i <
  \omega} \sigma_{a}(i)/4^{k_{i}}$ and $u_{b} = \sum_{i <
  \omega} \tau_{b}(i)/4^{l_{i}}$. By multiplying both sides of (\ref{PO}) with an appropriate integer, we can assume $p_{a}, q_{b} \in \egesz$ $(a < m,~b < n)$.  Let $j < \omega$ satisfy  \Keq\label{POI} 3 \cdot \max\b\{\sum_{a < m} |p_{a}|, \sum_{b < n} |q_{b}|\j\} < j.\Zeq It is enough to show that for every $j < h < \omega$, $\sum_{a <  m} p_{a}\sigma_{a}(h) = 0$; then the sums in (\ref{PO}) have a rational value. 

So suppose that for some $j < h < \omega$, $\sum_{a <  m} p_{a}\sigma_{a}(h)  \neq 0$. We consider every real in its base 4 decimal expansion, and for every $c < \omega$, the $c^{\textrm{th}}$ digit of $r \in \real$ is the coefficient of $4^{-c}$ in this expansion of $r$.

 By (\ref{POI}) and $j <h$, we have $3 \cdot\sum_{a < m} |p_{a}| < 4^{k_{h} - l_{h-1}}$, so there are $l_{h-1} < c_{0} < c_{1} \leq k_{h}$ such that the $c_{0}^{\textrm{th}}$ and $c_{1}^{\textrm{th}}$ digits of $\sum_{a < m} p_{a}f_{a}=\sum_{i <  \omega} \b(\sum_{a < m} p_{a}\sigma_{a}(i)\j)/4^{k_{i}} $ in its base 4 expansion differ. Also by $j <h$, we have $ 3 \cdot\sum_{b < n} |q_{b}| < 4^{l_{h} - k_{h}}$. So the $d^{\textrm{th}}$ digits of $$\sum_{j < n} q_{b}u_{b} = \sum_{i <  \omega} \b(\sum_{b < n} q_{b}\tau_{b}(i)\j)/4^{l_{i}}$$ for $l_{h-1} < d \leq k_{h}$ are either all 0 or all 3. This contradicts (\ref{PO}), so the proof of $\<F\>_{\rac} \cap \<U\>_{\rac} \subs \rac$ is complete.

Now $W \cup \{w_{n} \colon n <\omega\}  \ss U' \ss U$ so $\<F\>_{\rac} \cap \<W \cup \{v_{n} \colon n <
\omega\}\>_{\rac} \ss \rac$. Also by $W \cup \{w_{n} \colon n <\omega\} \ss U'$  we have $\<W \cup \{v_{n} \colon n <
\omega\}\>_{\rac} \cap \rac = \{0\}$. To summarize,  $\<F\>_{\rac} \cap \<W \cup \{v_{n} \colon n <
\omega\}\>_{\rac}  = \{0\}$, as stated.
\ep

\medskip

Once we have the compact set, its translates will be coded by the members of an almost disjoint family in $[\omega]^{\omega}$ of size $\kappa$. In the end we will need the following amended version of Solovay's Lemma.

  \bl\label{Solo}$(MA_{\kappa}({\sigma}\textrm{-centered}))$ Let
$\mc{A} \ss  [\omega]^{\omega}$ be an almost disjoint family of size $\kappa$.
 Let $\mc{B} \ss \mc{A}$ and suppose that for every $A \in\mc{B}$ a set $C_{A} \in [A]^{\omega}$ is given.  Then there exists $X \in
    [\omega]^{\omega}$ such that
    \begin{enumerate}
     \item $\max(X\cap A)\in C_A$ for $A\in \mc{B}$;
    \item $|X\cap A| = \omega$ for $A\in \mc{A} \sm \mc{B}$.
    \end{enumerate}
  \el
\prf Let
\begin{displaymath}
   P=\{\<x,b\> \colon x\in [\omega]^{<\omega},\, b\in
   [\mc{B}]^{<\omega},~ \max(x\cap B)\in C_B \text{ for } B\in b \},
  \end{displaymath}
and put $\<x,b\>\leq_{P} \<x',b'\>$ if and only if
$x' \ss x$, $b' \ss b$ and $x\cap B'= x'\cap B'$ for each $B'\in b'$.
 Since the  conditions
$\<x,b_0\>, \<x,b_1\>, \dots ,\< x,b_{n-1}\>$ have the joint extension
$\<x, b_0\cup b_1\cup\dots \cup b_{n-1}\>$, $P = \bigcup\{\{\<x,b\>\colon b \in [\mc{B}]^{< \omega}\} \colon x \in [\omega]^{<\omega}\}$ shows that $\<P,\le_P\>$ is ${\sigma}$-centered.

For every $B\in \mc{B}$ the set  $D_B=\{\<x,b\> \colon B \in b\}$ is dense in
$P$ since if $B\notin b$ then we have $n\in C_B\sm \max(B\cap \cup b)$ and $\<x\cup \{n\}, b\cup\{B\}\>\le \<x,b\>$ is in $D_B$.

For every $A\in \mc{A}\sm \mc{B}$ and $m < {\omega}$ the set
$D_{A,m}=\{\<x,b\> \colon \max(x\cap A)\ge m\}$ is dense in $P$ since for $n\in (A\sm m)\sm \cup b$, $\<x\cup\{n\},b\>\le
\<x,b\>$ is in $D_{A,m}$.

If $G$ is a $\{D_B \colon B\in\mc{B}\}\cup\{D_{A,m} \colon A\in \mc{A}\sm \mc{B},  m <
{\omega}\}$-generic filter then $X=\bigcup\{x \colon \<x,b\>\in G\}$ satisfies
the requirements.\ep

\medskip

\textbf{Proof of Theorem \ref{MA}.} 
 We can apply Lemma \ref{sov} to get $(k_{i})_{i < \omega} \ss \omega \sm \{0\}$, $(v_{n})_{n < \omega} \ss [0,1]$, $W\ss [0,1]$ and  $ F = \b\{ \sum_{i < \omega} 
\sigma(i)/4^{k_{i}} \colon \sigma \in \Sigma \j\}$ such that 
(\ref{ind00}) and (\ref{ind0}) from  Lemma \ref{sov} hold.
Let $\mc{A} \ss  [\omega]^{\omega}$ be an almost disjoint family of size
$\kappa$ and for every $A \in \mc{A}$ set $x(A) = \sum_{i < \omega}
\chi_{A}(i)/4^{k_{i}}$. Recall $\Sigma = 4^{\omega}$, and for every $n <
\omega$ let $\Sigma_{n} = \{\sigma \in \Sigma \colon n = \max\{ i < \omega
\colon \sigma(i) = 2\}\}$.  

For every $n < \omega$ let $F_{n} = \b\{ \sum_{i < \omega} \sigma(i)/4^{k_{i}}
\colon \sigma \in \Sigma_{n} \j\}$ and set $$K = W \cup F \cup
\bigcup\{F_{n} + v_{n} \colon n < \omega\}.$$ Note that $F_{n}$ $(n < \omega)$ are closed, \beq\label{ind}F \cup
\bigcup\{F_{n} + v_{n} \colon n < \omega\} \ss [0,2],~ 0 \in K \ss [0,2],\eeq and $\lim _{n
  \rar \infty} v_{n} = 0$ implies $\lim _{n
  \rar \infty} F_{n} + v_{n} = F$, hence $K$ is a compact set.
 We define
$$K_{n,A} =K + x(A)-v_{n}~(A \in \mc{A},~n < \omega)$$ and $\gotH_{0} = \{
K_{n,A} \colon A \in \mc{A},~n < \omega\}$. 

Set $Z=\{z \in F \colon
|\gotH_0(z)| = \kappa\}$. Let $\gotH_{1}$ consist of all translates of $K$ which avoid $Z$ and do not show up in $\gotH_{0}$, i.e.\ $$\gotH_{1}= \{K+d \colon (K+d) \cap  Z = \es, ~d \neq x(A)-v_{n}~(A \in \mc{A},~n < \omega)\}.$$ We show that the simple closed cover $\gotH = \gotH_{0}\cup  \gotH_{1}$  of $\real$ is $\kappa$-fold and has no two disjoint subcovers over $F$.

Pick an arbitrary $x \in \real$. If $|\gotH_0(x)| \geq \kappa$ we are done; so suppose $|\gotH_0(x)| < \kappa$. If for every $w \in W$, $(K+x-w)\cap Z = \es$ then by the definition of $\gotH_{1}$, $|\gotH_{1}(x)| \geq \kappa$. Similarly, if for every $f \in F$, $(K+x-f)\cap Z = \es$ then again $ |\gotH_{1}(x)| \geq \kappa$. If these cases fail to happen, then there are $w \in W$, $f \in F$, $y_{1}, y_{2} \in K$ and $z_{1}, z_{2}\in Z$  such that $z_{1} = y_{1}+x-w$ and $z_{2} = y_{2}+x-f$. Thus $x = z_{1}+w - y_{1}= z_{2}+f-y_{2}$. 

Let $D$ be a Hamel basis of $\real$ extending $W \cup \{v_{n} \colon n < \omega\}$. By (\ref{ind0}) of Lemma \ref{sov}, for every $y \in K$, if $y \in W$ then the expression of $y$ in the Hamel basis $D$ is $y$, while if $y \notin W$ then in the expression of $y$ in the Hamel basis $D$ no member of $W$ appears.

Consider the expression of $x$ in the Hamel basis $D$. We have $y_{1}, y_{2} \in K$ and $z_{1}, z_{2}, f \in F$, so in particular, by (\ref{ind0}) of Lemma \ref{sov} we have $z_{1}, z_{2}, f \in K \sm W$. Thus in the expression of $z_{2}+f-y_{2}$ in the Hamel basis $D$ no member of $W$ appears with positive coefficient. This implies $y_{1} = w$ hence $x=z_{1}$. Thus $x \in Z$ and so $|\gotH_0(x)| = \kappa$.

It remains  to see that $\gotH$ has no two disjoint subcovers over $F$. Let $c \colon \gotH_{0}\rar 2$. We find an $\ve \in \{0,1\}$ and an $x \in F$ such that $|\gotH_0(x)| = \kappa$, $\gotH_1(x) = \es$ and for every $A \in \mc{A}$ and $n < \omega$, $x \in K_{n,A}$ implies $c(K_{n,A}) = \ve$. This will complete the proof.

For each $A\in \mc{A}$ there exists  an ${\varepsilon}_A \in \{0,1\}$ and a $C_A\in [A]^{\omega}$ such that
$c(K_{n,A})={\varepsilon}_A$ for $n \in C_A$. Then there is ${\varepsilon} \in \{0,1\}$ and
$\mc{B}\in [\mc{A}]^{\kappa}$ such that ${\varepsilon}_B={\varepsilon}$ for $B\in \mc{B}$.

By applying Lemma \ref{Solo} we obtain $X \in [\omega]^{\omega}$ satisfying $\max(X\cap A)\in C_A$ $(A\in \mc{B})$ and $|X\cap A| = \omega$ $(A\in \mc{A} \sm \mc{B})$. Let $x=\sum_{i < \omega} (1+2 \chi_{X}(i)) / 4^{k_{i}}$, i.e.\ $x \in F$  and $x$ has digits 1 and 3 only. We show that this $x$ fulfills the requirements.

First we show $|\gotH_0(x)| = \kappa$. For every $z \in \real$ and $j < \omega$, let $z[j]$ denote the coefficient of $4^{-j}$ in the base $4$ expansion of $z$. For every $A \in \mc{A}$ and $j < \omega$ we have \Keq\label{AR}[x - x(A)][j]=\b\{ \begin{array}{ll} 3,&\textrm{if } j = k_{i} \textrm{ with } i \in X \sm A; \\ 2,&\textrm{if } j = k_{i} \textrm{ with } i \in X \cap A; \\ 1,&\textrm{if } j = k_{i} \textrm{ with } i \in \omega \sm (X \cup A); \\ 0,&\textrm{if } j = k_{i} \textrm{ with } i \in A \sm X. \end{array}\j. \Zeq  Thus for each $A\in \mc{B}$, $x-x(A) \in F_{\max (X \cap A)}$, hence $$x \in F_{\max (X \cap A)} +v_{\max (X \cap A)}+ x(A) -v_{\max (X \cap A)} \ss K_{\max (X \cap A), A}$$ and so $|\gotH_0(x)| = \kappa$. In particular,  $x \in Z$ and so $\gotH_1(x) = \es$.

It remains the show that for every $A \in \mc{A}$ and $n < \omega$, $x \in K_{n,A}$ implies $c(K_{n,A}) = \ve$.  Suppose $x \in K_{n,A}$ for some $A \in \mc{A}$ and $n < \omega$, i.e.\ \begin{multline}\notag x \in K+x(A)-v_{n} =  (W+x(A)-v_{n}) \cup \\ (F +x(A)-v_{n}) \cup
\bigcup\{F_{m} + v_{m}+x(A)-v_{n} \colon m < \omega\}.\end{multline}  By $x, x(A) \in F$ and (\ref{ind0}) of Lemma \ref{sov}, $$x \notin W+x(A)-v_{n},~x \notin F+x(A)-v_{n} ,~x \notin F_{m}+v_{m}+x(A) - v_{n} ~(m \neq n)$$ hence $x \in F_{n} +x(A)$. By (\ref{AR}), for $A\in \mc{A} \sm \mc{B}$ we have $x-x(A) \notin \bigcup _{n < \omega} F_{n}$ . Thus $x \in K_{n,A}$ implies $A \in \mc{B}$. Again by (\ref{AR}) we have  $n = \max(X \cap A) \in C_{A}$ so $c(K_{n,A}) = \ve$. This completes the proof in $\real$.

If $X = \omega^{\omega}$ or $X = 2^{\omega}$ take a continuous surjective map $\vp \colon X \rar [0,2]$ and set $\gotH_{X} = \{\vp^{-1}(H) \colon H \in \gotH\}$. Then $\gotH_{X}$ is clearly a $\kappa$-fold closed cover of $X$. It is a simple cover since by $K \ss [0,2]$, if $d_{1}, d_{2} \in \real$ satisfy $d_{1} \neq d_{2}$ and $(K+d_{1}) \cap [0,2] \neq \es$, $(K+d_{2}) \cap [0,2] \neq \es$ then $(K+d_{1}) \cap [0,2] \neq (K+d_{2}) \cap [0,2]$.  Since $\gotH$ has no two disjoint subcovers over $F$ and $F \ss [0,2]$, $\gotH_{X}$ has no two disjoint subcovers.\ep

\medskip

Corollary \ref{Cs_omega} implies in particular that in a positive partition result for closed covers the points covered   only by $
\omega$ many sets must be ignored.

\subsection{The Cohen real model}

In this section we will prove that in the Cohen real model every
 closed cover of the reals has an $[\omega_{1},\infty)$-good coloring. Note
 that by Corollary \ref{Cs_omega} it is impossible to get an 
${\omega}$-good
 coloring. Thus we have, in a sense, a best possible decomposition
 result. The proof is based on the weak Freeze-Nation property (see
 Proposition \ref{have} below), for which we need standard additional
 assumptions, such as GCH and $\doboz_{\lambda}$ for cardinals
 $\lambda$ with $ \cf(\lambda) = \omega$. 

Following \cite{Ku2}, we recall some notation. Let $V$ be our ground
model and let $\kappa$ be a cardinal.  We denote by $V^{C_{\kappa}}$
the model obtained from $V$ by adding $\kappa$ many Cohen reals  the
usual way. 

We will prove the following theorem.

\bt\label{Cs_Cohen} Suppose that GCH holds in $V$ and let $\kappa$ be
a cardinal. Suppose also that in $V$ we have  $\doboz_{\lambda}$  for
every cardinal $\lambda$ satisfying $\omega < \lambda \leq |\kappa|$,
$\cf(\lambda)=\omega$. In $V^{C_{\kappa}}$, let $(X, \tau)$ be a topological space which has a countable base, and let  $\gotH$ be a cover of $X$ by closed
sets. Then in $V^{C_{\kappa}}$  there exists an 
$[\omega_{1},\infty)$-good
coloring of $\gotH$. 
\et

The proof of Theorem \ref{Cs_Cohen} is based on the fact that in
$V^{C_{\kappa}}$ the poset $(\mc{P}(\omega), \ss)$ has the \emph{weak
  Freese-Nation property}. We recall it in  the following proposition
and we introduce the corresponding notion of good coloring on
$\mc{P}(\omega)$. 

\bc\label{have}(\cite[Theorem 15]{FS})
 Under the assumptions of Theorem \ref{Cs_Cohen}, in $V^{C_{\kappa}}$
 the poset  $(\mc{P}(\omega), \ss)$ has the weak Freese-Nation property,
 i.e.\ there is a 
  function $f \colon \mc{P}(\omega) \rar [\mc{P}(\omega)]^{\leq \omega}$ such
 that for every 
$A,B \in \mc{P}(\omega)$ with
$A \ss B$  there exists $C\in f(A)\cap f(B)$ satisfying $A\ss C\ss B$.
\ec

\bd\rm \label{Cs_max_wfn} Let $\mc{A}, \mc{B} \ss \mc{P}(\omega)$ be
arbitrary. A \emph{$\mc{B}$-good coloring} of $\mc{A}$
is a function $c \colon \mc{A} \rar \gotO$  such that for every $B \in
\mc{B}$, 
$|\mc{P}(B) \cap \mc{A}| \ge \omega_{1}$ implies
  $|\mc{P}(B) \cap \mc{A}| \ss c[\mc{P}(B) \cap \mc{A}]$. 
\ed

To get Theorem \ref{Cs_Cohen}, it is enough to prove the following theorem.

\bt\label{wfnthm}
Let $\mc{A}, \mc{B} \ss \mc{P}(\omega)$ be arbitrary. Assume
$(\mc{P}(\omega),\subseteq)$ has the weak Freese-Nation property. Then
$\mc{A}$ has a $\mc{B}$-good coloring. 
\et

\textbf{Proof of Theorem \ref{Cs_Cohen}.} By Proposition
\ref{mred} we can assume that $\gotH$ is simple, 
i.e.\ $\gotH=\mc{H}$.
Let $\{U_n:n<{\omega}\}$ be a base of $X$. For every closed set $Z \ss X$, define
  $B(Z)=\{n < \omega \colon U_n\cap Z=\es \}$. Since $Z\ss Z'$ if and only
 if $B(Z)\supseteq B(Z')$, $B$ is injective. 

Let $\mc{A}= \{B(H) \colon H\in\mc{H}\}$, 
$\mc{B} = \{B(\{x\}) \colon x \in
X\}$. By Theorem \ref{wfnthm} we have a $\mc{B}$-good coloring
$c^{\star} \colon \mc{A} \rar \gotO$. 
We show that $c \colon \mc{H} \rar \gotO$, $c=c^{\star} \circ B$  is 
an $[\omega_{1},\infty)$-good coloring of $\mc{H}$.

To see this, let $x\in X$ satisfy
$|\mc{H}(x)|\geq {\omega}_1$. Clearly, $B$ is a bijection between
$\mc{H}(x)$ and $\mc{P}(B(\{x\}))\cap \mc{A}$. Hence $|\mc{P}(B(\{x\}))\cap \mc{A}|\ge
{\omega}_1$ and so 
$c[\mc{H}(x)]=c^{\star}[\mc{P}(B(\{x\}))\cap \mc{A}] \supseteq
|\mc{P}(B(\{x\}))\cap \mc{A}|=
|\mc{H}(\{x\})|$, as required.\ep

\medskip

It remains to show Theorem \ref{wfnthm}.

 \medskip

\textbf{Proof of Theorem \ref{wfnthm}.} We prove  the statement by
induction on $\lambda = |\mc{A} \cup \mc{B}|$. 

If $\lambda \leq \omega$ an arbitrary coloring $c \colon \mc{A} \rar
\gotO$ works. 
Consider now $\lambda = \omega_{1}$. Enumerate $\mc{B}$ as
$\{B_{\alpha}  \colon \alpha < \omega_{1}\}$ such that each $B \in
\mc{B}$ occurs $\omega_{1}$ many times. We define $c \colon \mc{A}
\rar \omega_{1}$ by transfinite induction of length $\omega_{1}$,
extending $c$ to at most one further member of $\mc{A}$ at each step,
as follows. For every $B \in \mc{B}$ let $I_{B} = \{\alpha  <
\omega_{1} \colon B_{\alpha} = B\}$. In the $\alpha^{th}$ step of the
coloring if $\alpha \in I_{B}$ and $| \mc{P}(B) \cap \mc{A}| = \omega_{1}$
pick one $A \in \mc{A}$ such that $c(A)$ is not defined yet  and $A
\in \mc{P}(B)$. Define $c(A)=\tp(\alpha \cap I_{B})$. This coloring clearly
fulfills the requirements. 

Assume now that $\lambda > \omega_{1}$ and the statement holds for every
$\omega_{1} \leq \lambda' <{\lambda}$. 
Let $\mc{A},\mc{B}\ss \mc{P}({\omega})$ with $|\mc{A}\cup\mc{B}|=
{\lambda}$. 
Let $f \colon \mc{P}(\omega)\rar  [ \mc{P}(\omega)]^{\leq \omega}$ be a
function witnessing the 
    weak Freeze-Notion property of $\mc{P}(\omega)$. By closing $\mc{B}$
    under $f$ we can assume that $\mc{B}$ is $f$-closed. 

Let $\<M_{\alpha}:{\omega}_1\le{\alpha}<{\lambda}\>$ be a continuous, 
increasing sequence
of models of a large enough fragment of ZFC such that
$\mc{A}, \mc{B}, f \in M_{{\omega}_1}$,
$M_{\alpha}\in M_{{\alpha}+1}$,
$\alpha\ss M_{\alpha}$ and
  $|M_{\alpha}|=|{\alpha}|$ $({\omega}_1\le\alpha < \lambda)$. 
Let $\mc{A}_{\alpha}=\mc{A}\cap (M_{\alpha+1}\sm M_{\alpha})$, 
$\mc{B}_{\alpha}=\mc{B}\cap M_{\alpha+1}$.

For every
${\omega}_1\le \alpha < \lambda$
we have  $|\mc{A}_{\alpha}\cup \mc{B}_{\alpha}|=|{\alpha}|$.
So by the inductive hypothesis 
there is 
a coloring $c'_{\alpha} \colon
\mc{A}_{\alpha}\to |{\alpha}|$ which is 
$\mc{B}_{\alpha}$-good for $\mc{A}_{\alpha}$. 
By Lemma \ref{l:ha}, there is a function  $h_{\alpha}:|{\alpha}|\to {\alpha}$ such that 
 $h_{\alpha}[|{\alpha}|]={\alpha}$ and
 $ {\kappa} \ss h_{\alpha}[{\kappa}]$ for every cardinal ${\omega}\le {\kappa}<|{\alpha}|$.
Let $c_{\alpha}=h_{\alpha}\circ c'_{\alpha}$.
Then 
for every $B \in
\mc{B_{\alpha}}$, 
\ben[(i)]
\item \label{ab:i} $|\mc{P}(B) \cap \mc{A}_{\alpha}| \ge \omega_{1}$ implies
  $|\mc{P}(B) \cap \mc{A}_{\alpha}| \ss c_{\alpha}[\mc{P}(B) \cap \mc{A}_{\alpha}]$,
\item\label{ab:ii} $|\mc{P}(B) \cap \mc{A}_{\alpha}|=|{\alpha}|$ implies
  ${\alpha}\ss c_{\alpha}[\mc{P}(B) \cap \mc{A}_{\alpha}]$. 
\een

Let $c=\bigcup\{c_{\alpha} \colon {\alpha}<{\lambda}\}$; the
definition makes sense since for $\alpha \neq \beta$ we have
$\mc{A}_{\alpha} \cap \mc{A}_{\beta} = \es $. We show that $c$ is a
$\mc{B}$-good coloring. 

Assume on the contrary that there is $B \in \mc{B}$ such that
$|\mc{P}(B)\cap \mc{A}|\ge {\omega}_1$ but  $|\mc{P}(B) \cap \mc{A}| \not \ss
c[\mc{P}(B) \cap \mc{A}]$. Let $\alpha <\lambda$ be minimal such that we
can have such a $B$ in 
$M_{\alpha+1}$ and let ${\mu}\leq {\lambda}$ be an {\em uncountable
  regular} cardinal 
 such that
$|\mc{P}(B)\cap \mc{A}|\geq {\mu}$ but $\mu \not \ss c[\mc{P}(B)\cap
 \mc{A}]$. We distinguish three cases. 

Suppose first $|(\mc{P}(B)\cap \mc{A})\sm M_{\alpha}|\geq \mu$ and
$\alpha \geq \mu$. Then by  ${\mu} \ss M_{\alpha+1}$ we have $|((\mc{P}(B)\cap \mc{A})\sm M_{\alpha})\cap
M_{\alpha+1}|\geq 
\mu$. Hence $|\mc{A}_{\alpha}\cap \mc{P}(B)| \geq {\mu}$ and so
$c[\mc{A}\cap \mc{P}(B)]\supseteq
c_{\alpha}[\mc{A}_{\alpha}\cap \mc{P}(B)]\supseteq {\mu}$, a contradiction.

 Suppose next $|(\mc{P}(B)\cap \mc{A}) \sm M_{\alpha}|\geq \mu$ but
 ${\alpha}< {\mu}$. 
Let ${\sigma}\in {\mu}\sm c[\mc{P}(B)\cap\mc{A}]$ and let
 ${\beta}=\max({\alpha},{\sigma}+1)<{\mu}$. 
Then $|(\mc{P}(B)\cap \mc{A}) \sm M_{\beta}|\ge {\mu}$
and so $\omega_{1} \leq \beta \ss M_{\beta}$ implies $|\mc{P}(B)\cap
 \mc{A}_{\beta}|=|{\beta}|$. 
Thus $\beta \ss c_{\beta} [\mc{P}(B)\cap \mc{A}_{\beta}]$ by
(\ref{ab:ii}) and so $\sigma
 \in c_{\beta} [\mc{P}(B)\cap \mc{A}_{\beta}]\ss 
c [\mc{P}(B)\cap \mc{A}]$, a contradiction.

Finally  suppose $|(\mc{P}(B)\cap \mc{A}) \sm M_{\alpha}|< {\mu}$. With
$\nu = |M_{\alpha}\cap f(B)\cap \mc{P}(B)|\leq \omega$ enumerate
$M_{\alpha}\cap f(B)\cap \mc{P}(B)$ as $\{B_{i} \colon i < \nu\}$. 
For each $A\in \mc{A} \cap \mc{P}(B)\cap M_{\alpha}$ there is $B'\in
f(B)\cap f(A)$ with $A\ss B'\ss B$. Since $M_{\alpha}$ is
$f$-closed, $A \in M_{\alpha}$ implies $f(A) \ss M_{\alpha}$. Thus we have our $B'\in M_{\alpha}$, i.e.\ $B'=B_{n(A)}$ for some
$n(A) < \nu$. Therefore $$\mc{A} \cap \mc{P}(B)\cap M_{\alpha} =
\textstyle \bigcup_{n < \nu} \{A\in \mc{A}\colon  A\in
M_{\alpha},~ A\ss B_n \}.$$ 
Since $|\mc{A}\cap \mc{P}(B)\cap M_{\alpha}|\ge {\mu}$ there is
$n<\nu$ such that
$|\{A\in \mc{A}: A\in  M_{\alpha},~ A\ss B_n \}|\ge {\mu}$.
Since $B_n\in M_{{\alpha}^*+1}$ for some ${\alpha}^*<{\alpha}$, by the
minimality of ${\alpha}$ we have that $|\mc{A}\cap \mc{P}(B_n)|\ge {\mu}$
implies $\mu \ss c [\mc{A}\cap \mc{P}(B_n)] $. But
$ c[\mc{A}\cap \mc{P}(B_n)] \ss c[\mc{A}\cap \mc{P}(B)] $, a contradiction. This
completes the proof.\ep 

To close this section, we prove the following Lindel\"of-like property mentioned in the introduction.

\begin{ccc}\label{Marci}  Let $\gotH$ be a closed cover of $\real$ such that $|\gotH| < \cov(\mc{M})$.  Then $\gotH$ has a countable subcover of $\real$. In particular, for $\omega < \kappa < \cov(\mc{M})$, every $\kappa$-fold closed cover of $\real$ of cardinality $\kappa$ has a $\kappa$-good coloring.
\end{ccc}
\bp Let $\mc{U}$ be the collection of those open sets $U\ss \real$ for which $\gotH$ has a countable subcover of $U$; i.e. $$\mc{U} = \{U \ss \real \colon U \textrm{ is open, } \exists \mc{C} \in [\iH]^{\leq \omega} ~(U \ss \textstyle \bigcup \mc{C})\}.$$ Then $V = \bigcup \mc{U}$ is open. Since $\real$ is  hereditarily Lindel\"of, there is $\mc{V} \in [ \mc{U}]^{\leq \omega}$ such that $V = \bigcup \mc{V}$; in particular,  $V \in \mc{U}$.

To complete the proof of the first statement, it is enough to show that $V = \real$. Suppose $V \neq \real$ and set $F = \real \sm V$. Then $F$ is a nonempty closed set and $\gotH$ is a cover of $F$. Since $|\gotH| < \cov(\mc{M})$, there is an $H \in \mc{H}$ such that $H \cap F$ is non-meager in $F$ in the relative topology on $F$. Thus there is an open set $U \ss \real$ such that $U \cap F\neq \es$ and $U \cap F \ss H \cap F$. To summarize, we obtained that $\gotH$ has a countable subcover of $V \cup U$. This contradicts the definition of $\mc{U}$.
 
The second statement immediately follows from the first statement, so the proof is complete.
\ep

\section{Convex sets in $\real^{n}$}\lab{s:convex}

\subsection{Arbitrary convex sets}\label{arbconv}

In this section we observe that Theorem \ref{MA} and Theorem \ref{Cs_Cohen} imply that it is independent of ZFC whether an uncountable-fold cover of $\real^{n}$ $(1< n < \omega)$ by isometric copies of one compact convex set can be split into two disjoint subcovers.

\bt\label{arbcont} Let $1 <n < \omega$. Under the assumptions of Theorem \ref{MA}, there exists a $\kappa$-fold simple closed cover of $\real^{n}$ by isometric copies of one compact convex set which cannot be decomposed into two disjoint subcovers.
\et
\bp By rescaling the construction for Theorem \ref{MA}, there is a compact set $K \ss [\pi/4,\pi/2]$ and a set of translations $T \ss [-\pi/4,\pi/4]$ such that $\gotK = \{K+t \colon t \in T\}$ is a $\kappa$-fold simple cover over $[\pi/4,\pi/2]$ which cannot be split into two subcovers over $[\pi/4,\pi/2]$.

Let $\mathbb{O} \in \real^{n-2}$ denote the origin. For every $t \in \real$ set
$$H(t) = \conv \{(\cos(\vt+t), \sin (\vt+t)) \colon \vt \in K\}  \times \{\mathbb{O}\}$$
and let $\gotH_{0}=\{H(t) \colon t \in T\}$. Set $Y = \{(\cos(\vt), \sin (\vt))\colon \vt \in [\pi/4,\pi/2]\} \times \{ \mathbb{O}\}$ and let $\gotH_{1}$ be a $\kappa$-fold simple cover of $\real^{n}\sm Y$ by isometric copies of $H(0)$ which do not intersect $Y$. Such a $\gotH_{1}$ clearly exists. Then $\gotH = \gotH_{0} \cup \gotH_{1}$ fulfills the requirements.
\ep

The consistency of the existence of $[\omega_{1}, \infty)$-good colorings for compact covers follows from Theorem \ref{Cs_Cohen}.

\subsection{Axis-parallel closed rectangles}\label{axpar}

\bt\label{apr} There exists a countable family $\mc{R}$ of
axis-parallel closed rectangles in $\real^{2}$ such that $\mc{R}$ is
an $\omega$-fold cover of  $\real^{2}$ without two disjoint
subcovers. 
\et

We prove Theorem \ref{apr} in two steps: first we find an
$\omega$-fold cover of an abstract space without two disjoint
subcovers, then we show how this cover can be realized using
axis-parallel closed rectangles in $\real^{2}$.

Let $X=({\omega+1})^{\omega}\cup ({\omega+1})^{<{\omega}}$. For each 
${\sigma}\in ({\omega+1})^{<{\omega}}$ and $n\le{\omega}$, set
$$C_{{\sigma}^\frown n}=\{{\sigma}\}\cup \{f\in ({\omega+1})^{\omega}\colon{\sigma}^\frown n\subs
f\},$$ and let $\mc{C}=\{C_{{\sigma}^\frown n} \colon {\sigma}\in 
({\omega+1})^{<{\omega}},~ n\le {\omega}\}$.

\bl \label{comb} $\mc{C}$ is an $\omega$-fold cover of $X$ which cannot be split into two disjoint subcovers.
\el

\prf Pick an  arbitrary $x \in X$. If $x \in ({\omega+1})^{\omega}$ 
then $x\in C_{x|_{ k}}$ for each $k>0$.
 If $x \in (\omega+1)^{<\omega}$ then $x \in C_{x ^{\frown} n}$ 
$(n < \omega)$ so $\mc{C}$ is an $\omega$-fold
cover of $X$, indeed.

Split $\mc{C} = \mc{C} _{0} \cup \mc{C}_{1}$ where $\mc{C} _{0} \cap
\mc{C}_{1} = \es$. We show that if $\mc{C}_{0}$ is a cover of $X$ then
$\mc{C}_{1}$ is not a cover of $X$. So suppose $X = \cup \mc{C}_{0}$.
We define inductively a sequence $s \in (\omega+1)^{\omega}$ such that
$C_{s|_{n+1}} \in \mc{C}_{0}$ $(n < \omega)$; then  
we get $\mc{C}_{1}(s) = \es$, which shows
that $\mc{C}_{1}$ is not a cover of $X$. Let $n < \omega$ and suppose that $s|_{ n}$ is defined.  
Since $s|_{ n} \in X$ and  $\mc{C}_0$ is a cover of $X$, we have $\mc{C}_0(s|_{
n})\ne \empt$. So there is an
$m \leq \omega$ for which $C_{(s|_{ n}) ^{\frown} m} \in
\mc{C}_{0}$. Defining $s(n) = m$ completes the inductive step and the
proof.\ep

\medskip

\textbf{Proof of Theorem \ref{apr}.} 
First we show that it is enough to  construct a bijection $\vp$ between $X$ and a 
closed subset $F$ of $\real^2$
such that for each $C\in \mc{C}$ there is an axis parallel
rectangle $\Phi(C)$ such that $\vp[C]=F\cap \Phi({C})$. Indeed, since $F$ is closed, $\real^2\setm F$ has a countable, ${\omega}$-fold 
cover $\mc{D}$ by  axis-parallel closed rectangles which are all disjoint from $F$.
 Then $\mc{D}\cup \{\Phi(C)\colon C\in \mc{C}\}$ is a countable  
${\omega}$-fold cover
 of $\real^2$ by axis-parallel closed recangles which by Lemma \ref{comb} cannot be split into two disjoint subcovers.

We construct  $F$ as $F_1\cup F_2$ such that 
$\vp[({\omega}+1)^{<{\omega}}]=F_1$ and 
$\vp[({\omega}+1)^{{\omega}}]=F_2$.
Let $F_1$ be a countable closed subset of  the closed line segment connecting the points $(-1,1), (0,2) \in \real^{2}$, 
and let $\vp|_{ ({\omega}+1)^{<{\omega}}}$ be an arbitrary
bijection between $({\omega}+1)^{<{\omega}}$ and $F_1$.

We need some preparation to construct $F_2$.
For each ${\sigma}\in ({\omega}+1)^{<\omega}$,
by induction on $|{\sigma}|$, we construct a closed interval
$I_{\sigma}\subs [0,1]$ as follows (see Figure 1).
We set $I_{\es} = [0,1]$.
If $I_{\sigma}$ is constructed then
we choose $I_{{\sigma}^\frown n}\subs I_{\sigma}$ for $n\le {\omega}$
such that 
\ben[{(I}1{)}]
    \item   $\max I_{\sigma ^{\frown} n} < \min
    I_{\sigma ^{\frown} n'}$  
 and  
$(\max I_{{\sigma}^\frown n}-
\min I_{{\sigma}^\frown n}) <2^{-|{\sigma}|}$ $(n< n' \leq \omega)$;

\item  $\lim_{n < \omega} \max I_{\sigma
    ^{\frown} n} = \min I_{\sigma
    ^{\frown}\omega}$, $\min  I_{\sigma ^{\frown}0} = \min I_{\sigma}$
	  and $\max I_{\sigma ^{\frown} \omega} = \max I_{\sigma}$. 
\een

 \bigskip

 \bigskip

\begin{center}

 \vspace{-30pt}

\includegraphics[width=11cm]{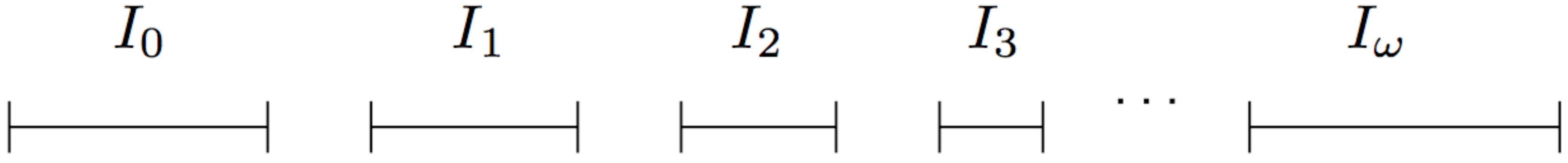}

{\small  $1^{\textrm{st}}$ level}

\vspace{5pt}

\includegraphics[width=13.5cm]{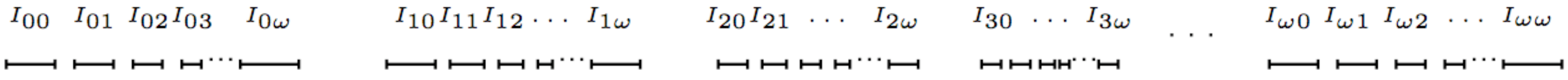}

{\small $2^{\textrm{nd}}$ level}
 \bigskip

Figure 1.
\end{center}

For every $x\in ({\omega}+1)^{\omega}$
set  $h(x) = \bigcap_{ l <{\omega}} I_{x|_{  l }}$. By (I2), for every $ l <{\omega}$ the set
$H_ l =\bigcup\{I_{\sigma}\colon {\sigma}\in ({\omega+1})^{l} \}$ is closed. So 
$H=\bigcap_{l <{\omega}} H_{ l } $ is also closed.

We define $\vp(x)=(h(x),h(x))$ for $x\in ({\omega}+1)^{{\omega}}$.
We have $H=h[({\omega}+1)^{{\omega}}]$, so $F_2=\vp[({\omega}+1)^{\omega}]=\{(z,z)\colon z\in H\}$ is closed in $\real^2$.

It remains to define $\Phi$. For 
every ${\sigma}\in  ({\omega}+1)^{<{\omega}}$ and $n\le {\omega}$
let 
$\Phi(C_{{\sigma}^\frown n})$ be the unique axis-parallel closed rectangle on $\real^{2}$ 
whose upper left corner is 
 $\vp({\sigma})$ and whose lower right corner is the point 
$(\max I_{{\sigma}\frown n},\min 
I_{{\sigma}\frown n})$.  Then for 
every ${\sigma}\in  ({\omega}+1)^{<{\omega}}$ and $n\le {\omega}$, $\Phi(C_{{\sigma}^\frown n})\cap F_1=\{\vp({\sigma})\}$
and  $$\Phi(C_{{\sigma}^\frown n})\cap F_2=\{(z,z) \colon 
z\in I_{{\sigma}\frown n}\}=\vp[\{x\colon {\sigma}^\frown n\subs x\}]= \vp[
C_{{\sigma}^\frown n}\cap ({\omega}+1)^{\omega}].$$
Thus $\vp[C_{{\sigma}^\frown n}]=\Phi(C_{{\sigma}^\frown n})\cap F$.
So the functions $\vp$ and $\Phi$ satisfy the requirements. \ep

\subsection{Polyhedra}\label{poly}

\newcommand{\opensets}{{\tau}}
\newcommand{\amal}[2]{\sqcup(#1\sqcap#2)}

The purpose of this section is to show that an uncountable-fold cover
of $\real^{n}$ by polyhedra has an $[\omega_{1},\infty)$-good coloring. We managed to obtain the following general result in this
direction, which allows us to treat covers by sets with very different
geometric constraints in a unified way. 

We introduce some notation in advance. Let $(X,\tau)$ be a topological
space, where $\tau$ stands for the family of all open subsets of $X$.
For $\gotB\subs \mc{P}(X)$ let
\begin{displaymath}
 \amal{\gotB}{\tau}=\{\cup\xcal \colon \xcal\in \br {B\cap G \colon B\in\gotB, G\in
  \opensets};{\omega};\} 
\end{displaymath}

\bt\label{genP} Let $(X,\tau)$ be a hereditarily Lindel\"of space
and let $\gotB\subs \mc{P}(X)$
be an intersection-closed family which is well-founded under $\ss$.
Then every cover  $\gotH$  with $\hcal\subs \amal \gotB{\tau}$
 has an $[\omega_{1},\infty)$-good coloring.
  \et

From Theorem \ref{genP} we have the following immediate corollaries.

\bcor\label{Pkov} Let $\kappa$ be an uncountable cardinal. Any $\kappa$-fold cover of $\real^{n}$
\ben
\item\label{Pkov1} by sets which can be obtained as countable unions of relatively open subsets of real affine varieties,
\item\label{Pkov2} by open or closed polyhedra,
\item\label{Pkov3} by open or closed balls,
\een  can be split into $\kappa$ many disjoint subcovers.
\ecor
\bp Since the polynomial ring of $n$ variables over the reals is
Noetherian, the family of real affine varieties in $\real^{n}$ is
intersection-closed and well-founded under $\ss$. So for \ref{Pkov1},
we can apply Theorem \ref{genP} with $\gotB$ standing for the real
affine varieties in $\real^{n}$. Statements \ref{Pkov2} and
\ref{Pkov3} are special cases of \ref{Pkov1}. 
\ep

We proceed to the proof of Theorem \ref{genP}.

\medskip

  \textbf{Proof of Theorem \ref{genP}.} By Proposition \ref{mred}, we can assume that $\gotH$ is a simple cover.
For each cardinal ${\lambda}$, let $(\circ_{\lambda})$ denote the
following statement:
\begin{itemize}
\item[$(\circ_{\lambda})$] If $\gotA=\{\acal_\alpha \colon \alpha
\in I\}\subs \br {B\cap G \colon B\in\gotB, G\in
  \opensets};{\omega};$, $|I|\le {\lambda}$ 
then there is a function $c \colon I\to
{\lambda}$ such that for each $x\in X$,
$|\{{\alpha} \colon x\in \cup\acal_{\alpha}\}|\ge {\omega}_1$ implies
$|\{{\alpha} \colon x\in \cup\acal_{\alpha}\}| \ss \{c({\alpha}) \colon x\in \cup\acal_{\alpha}\} $.  \end{itemize}
Thus it is enough to prove that $(\circ_{\lambda})$ holds for each ${\lambda}$.  
We do it by induction on ${\lambda}$.
Clearly, we can assume $I={\lambda}$.
For $\lambda \leq \omega$ any coloring fulfills the requirements. So let first ${\lambda}={\omega}_1$.

Take an arbitrary $\gotA = \{\acal_\alpha \colon \alpha   <
\omega_{1}\}$ and let $$\mc{B} = \{B \in \gotB \colon \exists \alpha <
\omega_{1}, ~\exists G \in \opensets ~(B \cap G \in \mc{A}_{\alpha})\}.$$
We have $|\mc{B}| = \omega_{1}$, so since $\gotB$ is well-founded, the
intersection-closed hull $\mc{B} ^{\cap}$ of $\mc{B}$ satisfies
$|\mc{B} ^{\cap}| = \omega_{1}$. Hence we can take an enumeration
$\mc{B} ^{\cap} = \{B_{\alpha} \colon  \alpha < \omega_{1}\}$. We also
fix a bijection $\vp \colon \omega_{1}\times \omega_{1} \rar
\omega_{1}$. 

For every $\al < \omega_{1}$ we construct inductively a countable
partial coloring $c_{\alpha} \colon {\omega}_1 \rar \gotO$, as
follows. Let $ \alpha < \omega_{1}$ and suppose that 
$c_{\eta}$ is defined for every $\eta < \alpha$. Set $\gotR_{\alpha} =
{\omega}_1 \sm \bigcup_{\beta < \alpha} \dom (c_{\beta})$.  Let
$\beta,\chi < \omega_{1}$ be such that $\vp(\beta,\chi) = \alpha$. For
$\gamma \in \gotR_{\alpha}$ let \Keq\label{Egy1} G(\gamma) =
\bigcup\{G \in \opensets \colon \exists B \in \gotB~(B_{\beta} \ss B,~B
\cap G \in \mc{A}_{\gamma})\},\Zeq and $$C_{\alpha} = \{x \in
B_{\beta} \colon |\{\gamma \in \gotR_{\alpha} \colon x \in G(\gamma)
\}| = \omega_{1}\}.$$ 

Since $X$ is hereditarily Lindel\"of, $C_{\alpha}$ is
Lindel\"of. Hence there is an $I_{\alpha} \in
[\gotR_{\alpha}]^{\omega}$ such that $C_{\alpha} \ss \bigcup_{\gamma
  \in I_\alpha} G(\gamma)$; thus by (\ref{Egy1}), $C_{\alpha} \ss
\bigcup \{\cup \mc{A}_{\gamma} \colon \gamma \in I_{\alpha}\}$, as
well. We define $c_{\alpha}$ by $\dom(c_{\alpha})=I_{\alpha}$ and 
$c_{\alpha}(\gamma) = \chi$ $(\gamma
\in I_{\alpha})$. This completes the $\alpha^{\textrm{th}}$ step of
the construction. We set $c = \bigcup_{ \alpha <
  \omega_{1} } c_{\alpha}$. 

We show
that $c$ witnesses $(\circ_{{\omega}_1})$. 
Suppose that $x \in X$ satisfies 
$|\{{\alpha} \colon x\in \cup\acal_{\alpha}\}|= {\omega}_1$, i.e.\
there are $I \in [\omega_{1}]^{\omega_{1}}$ and $B_{\alpha} \in
\gotB$, $G_{\alpha} \in \opensets$ $(\alpha \in I)$ such that $x \in
B_{\alpha} \cap G_{\alpha} \in \mc{A}_{\alpha}$ $(\alpha \in I)$. Let
$B = \bigcap_{\alpha \in I} B_{\alpha}$, then $B \in \mc{B}^{\cap}$,
that is $B = B_{\beta}$ for some $\beta < \omega_{1}$. Pick an
arbitrary $\chi < \omega_{1}$ and let $\alpha = \vp(\beta, \chi)$. 

Recall the construction of $c_{\alpha}$:  since $c_{\eta}$ $(\eta
<\alpha)$ are countable, we have $|I \cap \gotR_{\alpha}| =
\omega_{1}$. Hence $x \in C_{\alpha}$ and so $x \in \cup
\mc{A}_{\gamma}$ for some $\gamma \in I_{\alpha}$. Since $ c_{\alpha}(\gamma) = \chi$ and $\chi< \omega_{1}$
was arbitrary, the proof of the $\lambda = \omega_{1}$ case is
complete.

Let now ${\lambda}>{\omega}_1$ and suppose the statement holds for
every $\omega_{1} \leq \kappa < \lambda$. 
Let $\mc{M}=\<M,\in\>$ be a large enough model of a large enough
fragment of ZFC.
Let $\<M_{\alpha}\colon
{\omega}_1\le{\alpha} < {\lambda}\>$ be a continuous, increasing  chain of 
elementary submodels of $\mc{M}$ such that 
$|M_{\alpha}|=|{\alpha}|$ , $M_{\alpha}\in M_{{\alpha}+1}$, 
and $(X,\tau), \opensets,\gotB,\gotA\in M_{{\omega}_1}$.
For every set $y\in M_{{\lambda}}$, let 
$\mathrm{rank}(y)=\min\{{\alpha}\colon y\in M_{\alpha+1}\}$.

For every ${\omega}_1\le{\alpha} < {\lambda}$, 
let $J_{\alpha}={\lambda}\cap (M_{\alpha+1}\setm M_{\alpha})
=
\{\eta < {\lambda}\colon \mathrm{rank}(\eta )={\alpha}\}$. Then $|J_{\alpha}|=|{\alpha}|$.
By the inductive hypothesis, there is a coloring
$c'_{\alpha}\colon J_{\alpha}\to |\alpha|$ witnessing $(\circ)_{|{\alpha}|}$ for $\{\acal_{\xi} \colon {\xi}\in J_{\alpha}\}$.
By Lemma \ref{l:ha}, there is a function  $h_{\alpha} \colon |{\alpha}|\to {\alpha}$ such that 
 $h_{\alpha}[|{\alpha}|]={\alpha}$ and
 $ {\kappa} \ss h_{\alpha}[{\kappa}]$ for every cardinal ${\omega_{1}}\le {\kappa}<|{\alpha}|$.
Let $c_{\alpha}=h_{\alpha}\circ c'_{\alpha}$.
Then 
 for every $x \in X$, with ${\kappa}=|\{\eta \in J_{\alpha}\colon x\in \cup\acal_\eta\}|$,
\begin{enumerate}[(c1)]
\item ${\omega}_1\le \kappa
\le|{\alpha}|$ implies  $\kappa \ss 
\{c_{\alpha}(\eta)\colon x\in \cup\acal_\eta\}$;
\item $ \kappa = |{\alpha}|$ implies ${\alpha} = \{c_{\alpha}(\eta)\colon x\in \cup\acal_\eta\}.$
\end{enumerate}
 Set $c=\bigcup_{\alpha<\lambda}c_{\alpha} $; we show $c$
witnesses $(\circ_{\lambda})$. 

To this end, let $x\in X$ such that 
${\kappa}=|\{{\alpha} \colon x\in \cup\acal_{\alpha}\}|\ge {\omega}_1$.
Take ${\nu}_{\xi} \in \gotO$, $B_{\xi}\in\gotB$ and
$G_{\xi}\in\opensets$ $(\xi < {\kappa})$ such that $(\nu_{\xi})_{\xi < {\kappa}}$ are pairwise different and
$x\in B_{\xi}\cap  G_{\xi}\in
\acal_{{\nu}_\xi}$ $(\xi < {\kappa})$. Let ${\rho}_{\xi}=\mathrm{rank}({\nu}_{\xi})$ $({\xi}<{\kappa})$; we can assume $(\rho_{\xi})_{\xi<\kappa}$
is an increasing sequence.
Let ${\rho}=\sup \{{\rho}_{\xi}\stackrel.+1\colon {\xi}<{\kappa}\}$.
We can also assume that if ${\rho}$ is a successor
ordinal $\rho = {\rho}'+1$ then ${\rho}_{\xi}={\rho}'$ $({\xi}<{\kappa})$.

If ${\rho}$ is successor then
$|\{\eta \in J_{{\rho}'}\colon x\in \cup\acal_\eta\}|={\kappa}$
so we are done because $c_{{\rho}'}$ satisfies (c1). From now on assume
${\rho}$ is a limit ordinal. By the well-foundedness of $\gotB$ 
there are $B\in \gotB$ and  $F\in \br {\kappa};<{\omega};$
such that
\begin{equation}
\textstyle
B=\bigcap_{{\xi}<{\kappa}} B_{\xi}=\bigcap_{\xi \in F} B_{\xi}.
\end{equation}

Let ${\sigma}=\mathrm{rank}(F)\ge \mathrm{rank}(B)$.
Since ${\rho}$ is a limit ordinal we have  ${\sigma}<{\rho}$.
Thus
$\left|\left\{{\xi}<{\kappa}\colon  {\sigma} < {\rho}_{\xi} \right\}\right|={\kappa}$
and so
\Keq\label{Egy2} |\{\eta \in {\lambda}\setm \textstyle
M_{\sigma} \colon x\in \cup\acal_\eta\}|=\kappa . \Zeq For every $\alpha < \lambda$ let
$
G_\al=\bigcup\{G \in \opensets \colon \exists B' \in \gotB ~(B \ss B',~B' \cap G \in \acal_\al\}.
$ 

We distinguish two cases. First suppose ${\kappa}\leq {\sigma}$. 
Set \Keq\label{Egy3x} C = \{y \in B
\colon  |\{\eta \in {\lambda}\setm M_{\sigma} \colon y\in G_\eta\}| \geq
{\kappa}\}.\Zeq 
Since $C$ is Lindel\"of, by (\ref{Egy3x}) \begin{multline}
  \label{Egy4x}\textrm{there is a sequence }(K^{\star}({\zeta}))_{\zeta
    < {\kappa}} \textrm{ of pairwise disjoint countable} \\ \textrm{subsets of }
	{\lambda}\setm \textstyle M_{\sigma} \textstyle \textrm{ such that
}\textstyle C \ss \bigcup_{\eta \in K^{\star}(\zeta)} G_{\eta} 
\textrm{ for each }
~\zeta
< {\kappa}.\end{multline} 
But ${\kappa} \leq \sigma$ implies $B,
{\kappa} \in M_{\sigma+1}$, therefore $C \in M_{\sigma+1}$, as well. So by
elementarity (\ref{Egy4x}) holds in $M_{\sigma+1}$, i.e.\ there is a
sequence $(K(\zeta))_{\zeta < \kappa}$ of  pairwise disjoint countable
subsets of $\lambda \cap (M_{\sigma+1} \sm M_{\sigma})=J_{\sigma}$
in $M_{\sigma+1}$
such that  $C \ss \bigcup_{\eta \in K(\zeta)} G_{\eta}$
$(\zeta < \kappa)$. For every $\zeta < \kappa $ we have $K(\zeta)\subs
M_{\sigma+1}$ 
because $K(\zeta)\in M_{\sigma+1}$ and $K(\zeta)$ is countable. 
So $K(\zeta)\subs J_{\sigma}$ $(\zeta < \kappa)$. Since $G_\eta \cap B \subs \cup \acal_\eta$ $(\eta < \lambda)$, we have
 $C \ss \bigcup_{\eta \in K(\zeta)}\cup \acal_\eta$ $(\zeta < \kappa)$. 
Thus $\{\cup\acal_{\eta} \colon {\eta}\in J_{\sigma}\}$
is a ${\kappa}$-fold cover of $C$.
Since $x\in C$ and $c_{\sigma}$ satisfies $(c1)$, 
${\kappa}\subs \{c_{\sigma}(\eta) \colon x\in \cup\acal_\eta\}$, as required.

Finally suppose ${\sigma}<{\kappa}$. Fix an arbitrary $\beta \in
\gotO$ satisfying ${\sigma}<{\beta}<{\kappa}$. 
Set \Keq\label{Egy3xx} C = \{y \in B
\colon  |\{\eta \in {\lambda}\setm \textstyle 
M_{\beta} \colon y\in G_\eta\}| \geq
\beta\}.\Zeq 
 Since $C$ is Lindel\"of, by (\ref{Egy3xx}) \begin{multline}
  \label{Egy4xx}\textrm{there is a sequence }(K^{\star}({\zeta}))_{\zeta
    < \beta} \textrm{ of pairwise disjoint countable} \\ \textrm{subsets of }
	{\lambda}\setm M_{\beta} \textstyle \textrm{ such that
}\textstyle C \ss \bigcup_{\eta \in K^{\star}(\zeta)} G_{\eta} 
\textrm{ for each }\zeta
< \beta.\end{multline} But $\sigma \leq \beta$ implies $B,
\beta \in M_{\beta+1}$, therefore $C \in M_{\beta+1}$, as well. So by
elementarity (\ref{Egy4xx}) holds in $M_{\beta+1}$, i.e. there is a
sequence $(K(\zeta))_{\zeta < \beta}$ of  pairwise disjoint countable
subsets of $\lambda \cap (M_{\beta+1} \sm M_{\beta})=J_{\beta}$
in $M_{\beta+1}$
such that  $C \ss \bigcup_{\eta \in K(\zeta)} G_{\eta}$
$(\zeta < \beta)$. For every $\zeta < \beta$ we have $K(\zeta)\subs
M_{\beta+1}$ 
because $K(\zeta)\in M_{\beta+1}$ and $K(\zeta)$ is countable. 
So $K(\zeta)\subs J_{\beta}$ $(\zeta < \beta)$. 
Since $G_\eta \cap B \subs \cup \acal_\eta$ $(\eta < \lambda)$, we have
 $C \ss \bigcup_{\eta \in K(\zeta)}\cup \acal_\eta$ $(\zeta <
	\beta)$. 
Thus $\{\cup\acal_{\eta} \colon {\eta}\in J_{\beta}\}$ is a 
 $|\beta|$-fold cover of $C$.

 By (\ref{Egy2}) and by elementarity we
have $x \in C$, hence  $|\{\eta \in J_{\beta} \colon x\in
\cup\acal_\eta\}|=|{\beta}|$. 
Since $c_{\beta}$ satisfies (c2),
${\beta}\subs \{ c_{\beta}(\eta )\colon x\in \cup\acal_\eta\}\subs
\{ c(\eta )\colon x\in \cup\acal_\eta\}$. Since
${\beta}<{\kappa}$ was arbitrary, the proof is complete.

\ep

We remark that in the proof of Theorem \ref{genP}, formally we used
only the assumption that 
the space $X$ is a hereditarily ${\omega}_1$-Lindel\"of
space. However, a space is 
hereditarily ${\omega}_1$-Lindel\"of if and only if  it is
hereditarily Lindel\"of.

\section{Open problems}
\lab{s:open}

It is a matter of fact that whenever we considered the splitting problem of $\kappa$-fold covers for infinite $\kappa$ either we could establish the existence of a $\kappa$-good coloring or we could construct a $\kappa$-fold cover which cannot be split into two disjoint subcovers. Nevertheless, we could not prove that for $\kappa$-fold covers the existence of a 2-good coloring is equivalent with the existence of a $\kappa$-good coloring.

\bq \label{Q:gen}
Let $X$ be a set, $\kappa$ be an infinite cardinal and let $\mc{F} \ss \mc{P}(X)$ be arbitrary. Suppose every $\kappa$-fold cover $\gotH$ of $X$ satisfying $\mc{H} \ss \mc{F}$ has a 2-good coloring. Is it true then that every $\kappa$-fold cover $\gotH$ of $X$ satisfying $\mc{H} \ss \mc{F}$ has a $\kappa$-good coloring, as well?
\eq

In Section \ref{s:graph} we did not consider the splitting problem for hypergraphs.

\bq\label{hyperg} Examine the splitting problem of finite-fold and  infinite-fold edge covers of hypergraphs.
\eq

It would be interesting to know more on the consistency strength of the splitting of closed covers. In particular, one could examine whether maximal coloring of closed covers is possible in other well-known extensions than just the Cohen model. A special case is the following.

\bq \label{Q:ext}
Let $\kappa$ be an uncountable cardinal. Is it true in a random real extension of a model with GCH that every  $\kappa$-fold closed cover of $\real$ can be split into two disjoint subcovers?
\eq

We have seen that both under CH and under $\omega_{1} < \cov(\mathcal{M})$, an $\omega_{1}$-fold closed cover  $\gotH$ of $\real$ with $|\gotH| = \omega_{1}$ has an $\omega_{1}$-good coloring. However, we could not obtain it as a ZFC result.

\bq Is it consistent with ZFC that there exists an $\omega_{1}$-fold closed cover  $\gotH$ of $\real$ such that $|\gotH| = \omega_{1}$ but $\gotH$ cannot be split into two disjoint subcovers?
\eq

As we mentioned in the introduction, there are numerous open problems concerning the splitting of finite-fold covers of $\real^{n}$ by sets with special geometric properties. The interested reader is referred to \cite{PT} for more details. Here we propose problems for $\omega$-fold covers only.

\bq Is it true that every $\omega$-fold cover of $\real^{2}$ by \emph{translates} of one compact convex set can always be decomposed into two disjoint subcovers?
\eq

\bq Is it true that every $\omega$-fold cover of $\real^{n}$
\ben
\item  by translates or homothets of the unit cube,
\item by translates of the unit ball
\een
can be decomposed into two disjoint subcovers?
\eq

\noindent
\textbf{Acknowledgment.} We are indebted to P\'eter Er\-d\H os,
Andr\'as Frank, J\'anos Ger\-lits, Ervin Gy\H ori, Andr\'as Haj\-nal, Istv\'an Ju\-h\'asz,
L\'aszl\'o Lo\-v\'asz, Gyula  Pap, D\"om\"ot\"or P\'alv\"olgyi, G\'abor S\'a\-gi and Zolt\'an Szent\-mik\-l\'ossy for helpful
discussions.

We would also like to thank the referee for his numerous helpful comments and persistent work.

Our research was partially supported by the OTKA Grants K 68262, K 61600, K 49786 and F43620. We also gratefully acknowledge the support of \"Oveges Project of
\includegraphics[height=11pt]{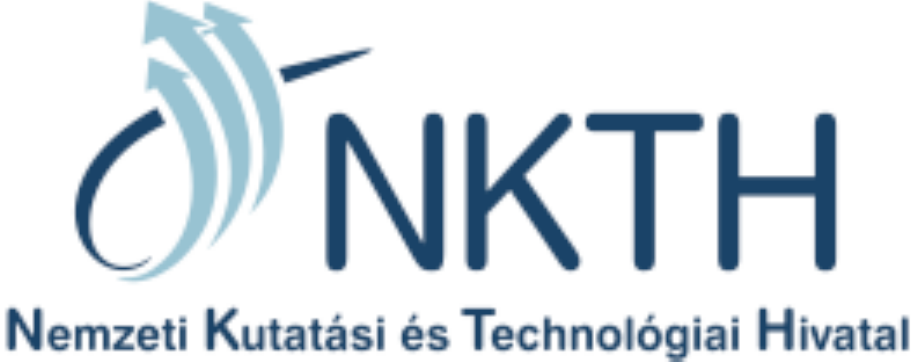} and
\includegraphics[height=11pt]{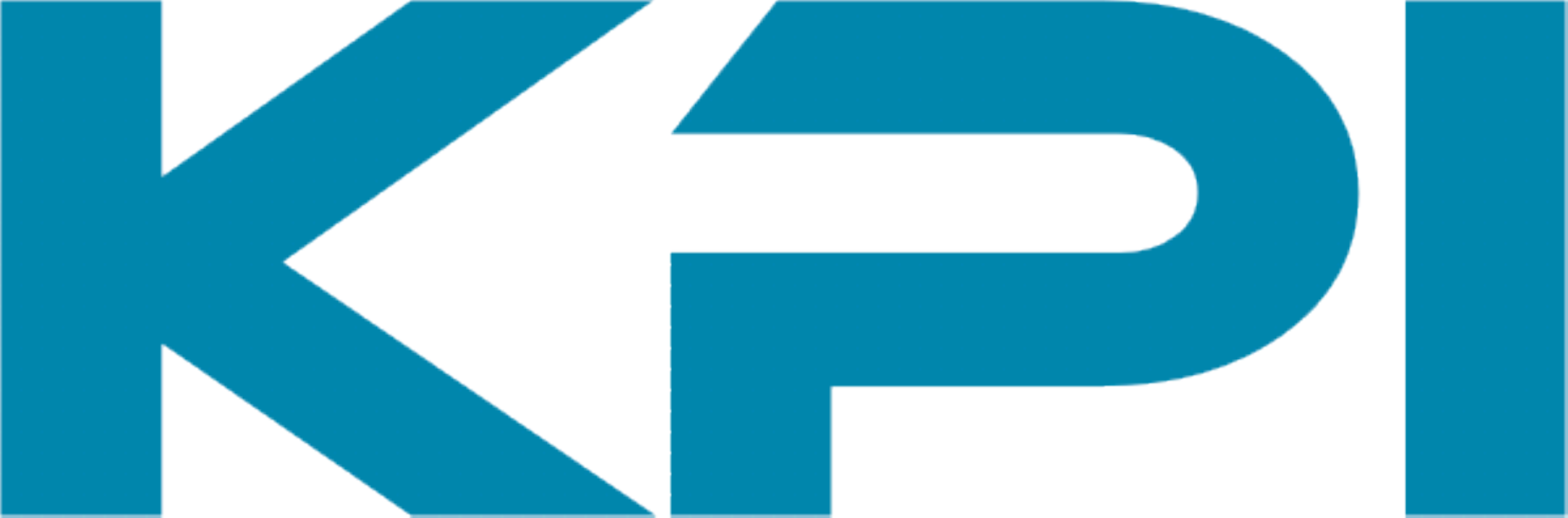}.

\bigskip

\noindent
\textsc{M\'arton Elekes}

\noindent
\textsc{Alfr\'ed R\'enyi Institute of Mathematics}

\noindent
\textsc{Hungarian Academy of Sciences}

\noindent
\textsc{P.O. Box 127, H-1364 Budapest, Hungary}

\noindent
\textit{Email:} \verb+emarci@renyi.hu+

\noindent
\textit{URL:} \verb+www.renyi.hu/~emarci+

\bigskip

\noindent
\textsc{Tam\'as M\'atrai}

\noindent
\textsc{Alfr\'ed R\'enyi Institute of Mathematics}

\noindent
\textsc{Hungarian Academy of Sciences}

\noindent
\textsc{P.O. Box 127, H-1364 Budapest, Hungary}

\noindent
\textit{Email:} \verb+matrait@renyi.hu+

\noindent
\textit{URL:} \verb+www.renyi.hu/~matrait+

\bigskip

\noindent
\textsc{Lajos Soukup}

\noindent
\textsc{Alfr\'ed R\'enyi  Institute of Mathematics}

\noindent
\textsc{Hungarian Academy of Sciences}

\noindent
\textsc{P.O. Box 127, H-1364 Budapest, Hungary}

\noindent
\textit{Email:} \verb+soukup@renyi.hu+

\noindent
\textit{URL:} \verb+www.renyi.hu/~soukup+

\end{document}